\documentclass[times,3p]{elsarticle}

\usepackage{amssymb}

\usepackage{lineno}

\usepackage{amsmath,amssymb,amsthm}
\usepackage{mathtools}
\usepackage[usenames,dvipsnames]{xcolor}
\usepackage{times}
\usepackage[utf8]{inputenc}
\usepackage{longtable}
\usepackage{textcomp}
\usepackage{mathtools}
\usepackage{dsfont}
\usepackage{xspace}
\usepackage{color}
\usepackage{algorithm}
\usepackage{algpseudocode}
\usepackage{pgfplots}
\usepackage{tikz}
\usepackage{subfigure}
\usepackage{wrapfig}
\usepackage{import}

\usepackage[hidelinks]{hyperref}
\def\equationautorefname~#1\null{(#1)\null}

\theoremstyle{definition}
\newtheorem{definition}{Definition}

\makeatletter
\DeclareRobustCommand\onedot{\futurelet\@let@token\@onedot}
\def\@onedot{\ifx\@let@token.\else.\null\fi\xspace}

\def\eg{\emph{e.g}\onedot} 
\def\ie{\emph{i.e}\onedot} 
\def\cf{\emph{cf}\onedot}

\def\etal{\emph{et al}\onedot}

\makeatother

\definecolor{maingreen}{rgb}{0.333333333333333, 1.0, 0.498039215686275}
\definecolor{myOrange}{RGB}{255, 169, 87 }
\definecolor{myGreen}{RGB}{180, 255, 162  } 
\definecolor{myBlue}{RGB}{100, 180, 255  } 
\definecolor{gnuplotbrown}{RGB}{255,76,0}
\definecolor{gnuplotblue}{RGB}{0,0,255}
\definecolor{gnuplotred}{RGB}{255,0,0}
\definecolor{gnuplotpurple}{RGB}{255,2,255}
\definecolor{gnuplotgray}{RGB}{131,128,126}
\definecolor{uniblau}{HTML}{004291}
\definecolor{bigsblau}{HTML}{365079}
\definecolor{bigsblau50}{HTML}{9CA7BC}
\definecolor{bigsblau25}{HTML}{CDD3DD}
\definecolor{uniorangedark}{HTML}{E6B400}
\definecolor{uniorange}{HTML}{FFCB0E}
\definecolor{uniorange!50}{HTML}{FFE586}
\definecolor{uniwhite}{HTML}{FFF2C2}
\definecolor{hcmgruen}{HTML}{567877}
\definecolor{hcmgruen50}{HTML}{AFBDBE}
\definecolor{hcmgruen25}{HTML}{D7DEDE}
\definecolor{himgrau}{HTML}{626566}
\definecolor{himgrau75}{HTML}{949592}
\definecolor{himgrau50}{HTML}{C5C4BE}
\definecolor{himgrau25}{HTML}{F5F5F5}
\definecolor{textgrau}{HTML}{000000}
\definecolor{black}{HTML}{000000}
\definecolor{white}{HTML}{FFFFFF}
\definecolor{lightgray}{gray}{0.9}
\definecolor{lila}{HTML}{FF34B3}
\definecolor{lightblue}{HTML}{99D6FF}
\colorlet{blue}{uniblau}
\colorlet{greyblue}{bigsblau}
\colorlet{hcmgelb}{uniorange}
\colorlet{yellow}{uniorange}
\colorlet{tafelgruen}{hcmgruen}
\colorlet{green}{hcmgruen}

\newcommand{\red}[1] {{\color{red}{{#1}}}}

\newcommand{\notinclude}[1]{}

\IfFileExists{todonotes.sty}{%
	\usepackage{todonotes}
	\usepackage{comment}
	\usepackage{setspace}
	\let\todoavailable1	
}{}%

\providecommand\MR{} 

\ifdefined\todoavailable
	\makeatletter
	\newcommand{\inlinetodo}[2][]{%
		\@todo[caption={}, inline, #1]{\begin{spacing}{0.5}#2\end{spacing}}%
	} 
	\newcommand{\missing}[1]{ \@todo[color=red!40]{Missing: #1} }
	 
	\newcommand{\BH}[2][]{\@todo[color=yellow, #1]{BH: #2}}
	\newcommand{\KH}[2][]{\@todo[color=yellow, #1]{KH: #2}}
	\newcommand{\JS}[2][]{\@todo[color=yellow, #1]{JS: #2}}
	\renewcommand{\MR}[2][]{\@todo[color=yellow, #1]{MR: #2}}
	\makeatother
\else
	
	\newcommand{\missing}[1]{ {\red{missing: [} #1 \red{]!!!}} }
	
	\newcommand{\BH}[2][]{\red{BH: [} #2 \red{].}}
	\newcommand{\KH}[2][]{\red{KH: [} #2 \red{].}}
	\renewcommand{\MR}[2][]{\red{MR: [} #2 \red{].}}
	\newcommand{\JS}[2][]{\red{JS: [} #2 \red{].}}
	\newcommand{\inlinetodo}[2][]{{\red{TODO: {#2} }\hfill\\}}%
	
\fi

\newcommand{\R}{\mathbb{R}}
\newcommand{\N}{\mathbb{N}}

\newcommand{\tr}{\mathrm{tr}\,}
\DeclareMathOperator*{\Id}{Id}
\renewcommand{\bmod}{{\,\mathrm{mod}\,}}

\DeclareMathOperator*{\argmin}{arg\,min}
\renewcommand{\d}{\,\mathrm{d}}
\newcommand{\SO}{\mathit{SO}}

\newcommand{\manifold}{\mathcal{M}}
\newcommand{\vertexSpace}{\mathcal{N}}
\newcommand{\imm}{X}
\newcommand{\immAlt}{\tilde{\imm}}
\newcommand{\immEx}{\imm^\ast}

\newcommand{\edge}{e}
\newcommand{\vertex}{v}
\newcommand{\face}{f}
\newcommand{\pos}{X}
\newcommand{\area}{a}
\newcommand{\normal}{N}
\newcommand{\vertices}{\mathcal{V}}
\newcommand{\edges}{\mathcal{E}}
\newcommand{\faces}{\mathcal{F}}
\newcommand{\numV}{|\vertices|}
\newcommand{\numE}{|\edges|}
\newcommand{\numF}{|\faces|}
\newcommand{\edgeSum}{\sum_{\edge \in \edges}}
\newcommand{\nv}{n_v}
\newcommand{\curve}{{\z}}

\newcommand{\graph}{\mathcal{G}}
\newcommand{\connec}{\graph}
\newcommand{\simplicial}{\mathcal{K}}
\newcommand{\underlyingSpace}{\vert \simplicial \vert}

\newcommand{\W}{\mathcal{W}}
\newcommand{\mem}{{\mbox{{\tiny mem}}}}
\newcommand{\bend}{{\mbox{{\tiny bend}}}}

\newcommand{\Wmem}{\W_\mem}
\newcommand{\Wbend}{\W_\bend}

\newcommand{\z}{{z}}
\newcommand{\len}{{l}}
\newcommand{\dih}{{\theta}}
\newcommand{\projZ}{{Z}}
\renewcommand{\frame}{{F}}

\newcommand{\intRot}{R}
\newcommand{\thickSqr}{\delta^2}
\newcommand{\metric}{g}
\newcommand{\Hess}{\mathrm{Hess}\,}

\newcommand{\dfFF}{G}

\newcommand{\dDisTen}{\mathcal{G}}

\newcommand{\trineq}{\mathcal{T}}
\newcommand{\trineqSet}{\mathcal{Z}_\mathcal{T}}
\newcommand{\dint}{\mathcal{I}}
\newcommand{\qint}{\mathcal{Q}}

\newcommand{\memw}{\alpha}
\newcommand{\benw}{\beta}

\newcommand{\rec}{\mathcal{R}}

\newcommand{\energy}{\mathcal{E}}
\newcommand{\lagrange}{L}
\newcommand{\mult}{{\lambda}}
\newcommand{\pen}{{\mu}}
\newcommand{\shift}{{\tau}}

\graphicspath{{./figures/}}

\journal{Computer Aided Geometric Design}

\begin{document}

\begin{frontmatter}



\title{Geometric optimization using nonlinear rotation-invariant coordinates}
\date{February 5, 2020}

\author[INS]{Josua Sassen}
\ead{josua.sassen@uni-bonn.de}

\author[INS]{Behrend Heeren}
\ead{heeren@ins.uni-bonn.de}

\author[Delft]{Klaus Hildebrandt}
\ead{K.A.Hildebrandt@tudelft.nl}

\author[INS]{Martin Rumpf}
\ead{martin.rumpf@ins.uni-bonn.de}

\address[INS]{Institute for Numerical Simulation,  Universität Bonn, Endenicher Allee 60, 53115 Bonn, Germany}
\address[Delft]{Delft University of Technology, EEMCS -- Dept. Intelligent Systems, Van Mourik Broekmanweg 6, 2628 XE Delft, Netherlands}

\begin{abstract}
Geometric optimization problems are at the core of many applications in geometry processing. 
The choice of a representation fitting an optimization problem can considerably simplify solving the problem. 
We consider the Nonlinear Rotation-Invariant Coordinates (NRIC) that represent the nodal positions of a discrete triangular surface with fixed combinatorics as a vector that stacks all edge lengths and dihedral angles of the mesh.
It is known that this representation associates a unique vector to an equivalence class of nodal positions that differ by a rigid body motion.
Moreover, integrability conditions that ensure the existence of nodal positions that match a given vector of edge lengths and dihedral angles have been established.
The goal of this paper is to develop the machinery needed to use the NRIC for solving geometric optimization problems. 
First, we use the integrability conditions to derive an implicit description of the space of discrete surfaces as a submanifold of an Euclidean space and a corresponding description of its tangent spaces.
Secondly, we reformulate the integrability conditions using quaternions and provide explicit formulas for their first and second derivatives facilitating the use of Hessians in NRIC-based optimization problems.
Lastly, we introduce a fast and robust algorithm that reconstructs nodal positions from almost integrable NRIC.
We demonstrate the benefits of this approach on a collection of geometric optimization problems. 
Comparisons to alternative approaches indicate that NRIC-based optimization is particularly effective for problems involving near-isometric deformations.   
\end{abstract}


%


\end{frontmatter}



\section{Introduction}
Geometric optimization problems are central to geometry processing as most methods involve some form of optimization as a step in their pipeline. 
Shape deformation problems are inherently nonlinear and therefore can be difficult to solve. In particular, problems involving near-isometric deformations are typically ill-conditioned due to the combination of high stretching and low bending resistance. 
Moreover, physical objectives are often invariant with respect to rigid body motions and the alignment of a mesh in Euclidean space is considered to be a post-processing task.
However,  this rigid body motion invariance can cause conceptual and numerical issues when working with nodal positions. 
Thus it is beneficial to find degrees of freedom for mesh description and corresponding deformations which are rigid body motion invariant. 

In this work, we study the Nonlinear Rotation-Invariant Coordinates (NRIC) that describe the immersion of a mesh using the edge lengths and dihedral angles of the mesh instead of the nodal positions.
Beyond their inherent invariance to rigid transformations, these coordinates offer additional benefits, such as their natural occurrence in discrete deformation energies and their representation of natural modes of deformation in a localized sparse fashion.
For example, when a human character (represented by a triangle mesh) lifts her straight arm, the induced variations in nodal positions comprise the entire arm. 
However, the same variation encoded in the change of lengths and angles is limited to the shoulder region, \ie the place where the actual physical work is done.  

Prior work on shape interpolation by \citet{WiDrAl10} and \citet{FrBo11} showed that linear blending of the NRIC for a set of shapes already yields interesting nonlinear deformations. 
However, since in general nodal positions that realize given edge lengths and dihedral angles may not exist, these methods rely on optimization in the space of nodal positions.  

\paragraph*{Contribution}
The basis of our approach is the triangle inequalities and the integrability conditions derived by \citet{WaLiTo12}. 
Our goal is to provide the machinery required to formulate and solve geometric optimization problems entirely in NRIC. 
\begin{itemize}
	\item We reformulate the integrability conditions using quaternions and use this to provide an implicit description of the NRIC manifold along with its tangent spaces.
	\item We reformulate the nonlinear energy from \cite{HeRuSc14} in NRIC and provide its derivatives to equip the NRIC manifold with a Riemannian metric.
	\item For solving (constrained) geometric optimization problems in the NRIC manifold, we describe an approach based on the augmented Lagrange method.
	In this context, we illustrate how to efficiently handle the triangle inequality constraints using the natural barrier term in the nonlinear energy and a modified line search.
	This also includes explicit formulas for the second derivatives of the integrability conditions, which are needed for evaluating the Hessian of objectives acting on NRIC.
	\item Finally, we introduce a hybrid algorithm to construct nodal positions of a discrete surface from NRIC which do not necessarily fulfill the integrability conditions. 
	The algorithm uses an adaptive mesh traversal algorithm as initialization to a Gauß--Newton solver. 
	In our experiments, this proves to effectively reduce the number of required Gauß--Newton iterations. 
	Typically, a single iteration is sufficient or even no iteration is needed.
\end{itemize}
Experiments demonstrate the utility of our framework for various applications such as geodesics in shape space and paper folding. 
Our approach is particularly well-suited to deal with near isometric deformations of discrete shell surfaces, which is underpinned by a variety of numerical examples.

\paragraph*{Organization} 
The remainder of this paper is organized as follows. 
After reviewing related work in Section \ref{sec:relatedWork}, we summarize the necessary background on the established discrete integrability conditions as introduced by \citet{WaLiTo12} in Section \ref{sec:background}. 
We define our NRIC manifold and reformulate the integrability conditions using quaternions in Section \ref{sec:nric}.
In \autoref{sec:energy}, we discuss the nonlinear deformation energy.
Afterwards, we introduce a corresponding variational calculus in Section \ref{sec:varProblems}. 
The robust reconstruction of nodal positions from lengths and angles is discussed in Section \ref{sec:recon}. 
Finally, we show a series of applications in Section \ref{sec:results} and discuss limitations and challenges in Section \ref{sec:discussion}.


\section{Related Work}\label{sec:relatedWork}
In this section, we discuss relevant work on linear and nonlinear coordinates, rigidity of triangle meshes, shape interpolation, shape spaces and near-isometric deformation. 

\paragraph*{Linear coordinates} 
For solving problems in geometry processing, it can be useful to switch from the usual nodal coordinates to a different representation that is adapted to the given task. We distinguish between coordinates that depend linearly and nonlinearly on the nodal coordinates. 
\emph{Differential coordinates} use discrete differential operators on a triangle mesh to define coordinates. Two examples are \emph{gradient-domain} approaches for meshes \cite{YuZhXuShBaGuSh04,SuPo04}, which operate on the gradients of functions, and the \emph{Laplace coordinates}~\cite{SoCoLiAlRoSe04,LiSoCoLeRoSe04}, which make use of the discrete Laplace--Beltrami operator. 
Since the differential coordinates depend linearly on the nodal positions, the immersion that best matches given differential coordinates can be found by solving a linear least-squares problem.  
While linearity of the coordinates facilitates computations, it also fundamentally limits their applicability. 
For example, shape editing approaches that use linear coordinates often yield unnatural and distorted shapes when larger deformations are involved~\cite{BoSo08}. 

\paragraph*{Nonlinear coordinates} 
In classical differential geometry, the fundamental theorem of surfaces~\cite{Do76} states that two immersion of a surface to $\mathbb{R}^3$ differ by a rigid motion if and only if the first and second fundamental forms agree and provides integrability conditions that guarantee the existence of an immersion for a given first and second fundamental form.  
This motivates using discrete analogs to the fundamental forms as coordinates for triangle mesh processing. 
Explicitly, the list of all edges lengths and dihedral angles is used. Analogous to the classical theorem, the nodal positions of two meshes with the same combinatorics agree up to a global rigid motion if and only if all edge lengths and all dihedral angles agree.
Integrability conditions that guarantee the existence of nodal positions realizing a given vector of edge lengths and dihedral angles were derived by \citet{WaLiTo12}. 
The integrability conditions are formulated using moving frames associated with the triangles of the mesh. 
Already in earlier work, \citet{LiSoLe05,LiCoGaLe07} used moving frames to define coordinates on triangle meshes. 
Our goal is to extend this line of work by providing the tools and structures needed for solving optimization problems that are formulated in the nonlinear coordinates.   

\paragraph*{Rigidity} 
While the existence and uniqueness results for the nonlinear coordinates require both, the edge lengths and the dihedral angles, rigidity results can already be obtained if only edge lengths are considered. 
For convex polytopes, Cauchy's and Dehn's rigidity theorems \cite{De16} show rigidity and infinitesimal rigidity and \citet{Gl75} showed that \emph{almost all} simply-connected polyhedra are rigid. 
An example of polyhedra that allow for isometric continuous deformations, which are non-rigid, is Cornelly's sphere \cite{Co77}. 
In this paper, we will formulate infinitesimal rigidity in terms of NRIC.
In recent work, \citet{AmRo18} studied the dihedral rigidity of polyhedra and parametrized triangle meshes via dihedral angles. 
Related to rigidity is the problem of computing an immersion from prescribed edge lengths. 
Algorithms for this problem were proposed by \citet{BoEyKoBr15} and 
\citet{ChKnPiSc18}. 
Relaxing the concept of rigidity, conformal geometry identifies metrics that differ only by a conformal factor. \citet{CrPiSc11} study the numerical treatment of the integrability conditions for surfaces in this setting. 

\paragraph*{Near-isometric deformations} 
Isometric and near-isometric deformations are important for computational folding of piecewise flat or developable structures \cite{KiFlChMiShPo08,BoVoGoWa16,StGrCr18,RaHoSo18b}. 
The computation of near-isometric deformation can be done by simulating elastic shells consisting of stiff material with low bending resistance \cite{BuZoGr06,SoVoWaGr12,NaPfBr13}. 
These materials yield ill-conditioned problems that are difficult to solve numerically. 
The NRIC perspective improves the numerical accessibility of such problems.

\paragraph*{Shape interpolation}
Shape interpolation, also called blending or morphing, is an important problem in geometry processing which is used for applications such as deformation transfer \cite{BaVlGrPo09,YaGaLaRoXi18}, motion processing \cite{PrKaChCoHo16}, example-based methods for shape editing \cite{FrBo11}, inverse kinematics \cite{SuPo04,Wa16}, and material design \cite{MaThGrGr11}. 
Approaches to shape interpolation based on linear coordinates use non-linear operations for blending the coordinates. 
For example, the gradient-domain approach of \citet{XuZhWaBa05} extracts the rotational components from deformation gradients via polar decomposition and applies nonlinear blending operations to these components. 
While the nonlinear blending helps to compensate for linearization artifacts, it is a difficult task to estimate the local rotations that resolve large deformations. 
\citet{KiGa08} and \citet{GaLaLiChShXi16} introduce improved nonlinear blending operations for the rotational components.
\citet{WiDrAl10} introduce a scheme for shape interpolation using nonlinear coordinates. Their method linearly blends edge lengths and dihedral angles 
and uses a multi-scale shape matching algorithm for constructing interpolating shapes. 
\citet{FrBo11} model the process of finding the shape that best matches the blended lengths and angles as a nonlinear least-squares optimization problem and solve it using a multi-resolution Gau\ss--Newton scheme. 
A related approach by \citet{WuBoShRoBr10} blends edge lengths and the normal vectors of two example shapes and constructs the intermediate shapes using a mesh traversal algorithm based on a minimal spanning tree with dihedral angle differences as weights.
Model reduction approaches that enable real-time shape interpolation have been introduced by \citet{TyScSe15} and \citet{RaEiSe16}.
Related to shape interpolation are shape spaces, which are shape manifolds equipped with a Riemannian metric. 
Shape spaces are used for various applications in computer vision, computational anatomy, and medical imaging. For a general introduction to shape space and their applications, we refer to the textbook of \citet{Yo10}.
\citet{KiMiPo07} introduced a Riemannian metric on spaces of triangle meshes and show that concepts from Riemannian geometry such as the exponential map and parallel transport can be used for geometry processing tasks like deformation transfer, shape interpolation, and extrapolation.  
\citet{HeRuWa12,HeRuSc14} propose an alternative physically-based metric on the shape space of triangle meshes that reflects the viscous dissipation required to physically deform a thin shell. 
\citet{BrTyHi16} derive a discrete curve shortening flow in shape space and use it for processing animations of deformable objects. 
While the shape spaces study deformations of meshes with fixed connectivity, functional correspondences \cite{OvBeSoBuGu12} can be used to blend \cite{KoBrBrGlKi13} and analyze \cite{RuOvAzBeChGu13} pairs of meshes with different connectivity. 
In recent work, the functional correspondences of intrinsic and extrinsic geometry \cite{CoSoBeGuOv17} and deformation fields \cite{CoOv19} have been studied.

\section{Background}\label{sec:background}
In this section, we briefly review the work by \citet{WaLiTo12} who introduced a discrete version of the fundamental theorem of surfaces. 
We consider a simplicial surface, which is a simplicial complex $\simplicial = (\vertices, \edges, \faces)$ consisting of sets of vertices $\vertices$, edges $\edges \subset \vertices \times \vertices$ and faces $\faces \subset \vertices \times \vertices \times \vertices$ such that the topological space $\underlyingSpace$ obtained by identifying each face with a standard two-simplex and gluing the faces along the common edges is a two-dimensional manifold. 
A map $\imm \colon \vertices \to \R^3$ is called \emph{generic} if for every face, the three vertices are in general position, \ie, there is no straight line in $\R^3$ containing the three vertices. 
We define 
\begin{align}\label{eq:vertexSpace}
\vertexSpace := \{ \imm(\vertices) \mid \imm \colon \vertices \to \R^3 \text{ generic } \}  \subset \R^{3\numV} \, \, ,
\end{align}
which we denote the space of \emph{discrete surfaces}. 
For any $\imm$, there is a unique map $\immEx\colon \underlyingSpace \to \R^3$ that is continuous, an affine map of each simplex, and interpolates $\imm$ at the vertices. 
$\immEx$ maps the faces of $\underlyingSpace$ to triangles in \(\R^3\), and, if $\immEx$ corresponds to a generic map $\imm$, none of the triangles degenerates. 
We will need this property to ensure that the elastic energies we consider in Section~\ref{sec:energy} are well-defined. 

\begin{wrapfigure}{r}{0.45\columnwidth}
	\centering
	\scalebox{0.8}{
		\Large 
		\def\svgwidth{90mm}
		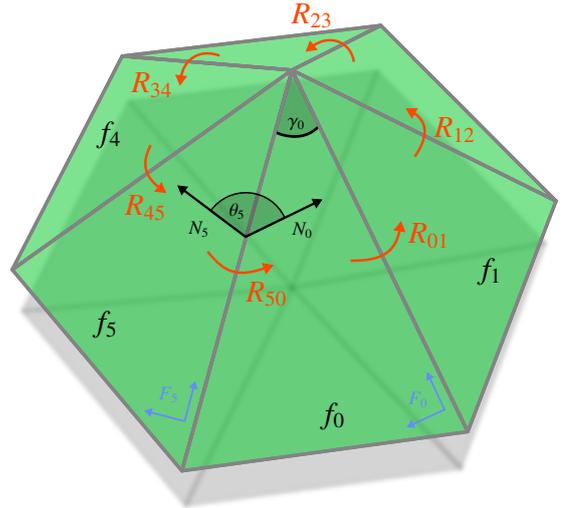
	}
	\caption{
		Construction of integrability condition for the 6-loop of faces around a vertex. 
			Specifically, the transition rotation \(\intRot_{50} = R_0(\theta_5) R_2(\gamma_{0})\) (orange) is constructed from the dihedral angle \(\theta_5\) and interior angle \(\gamma_0\) (both black).
			It transforms the frame \(\frame_5\) into frame \(\frame_0\) (both blue, normal vector not shown), \ie \(\frame_0 =  \frame_5 \intRot_{50}\).
			The other transition rotations are constructed in the same way and applying them sequentially yields the integrability condition \(\intRot_{01}\intRot_{12}\intRot_{23}\intRot_{34}\intRot_{45}\intRot_{50} \protect\overset{!}{=} \Id\).
	}
\end{wrapfigure}

Since we assume that the underlying simplicial complex remains unchanged, by abuse of notation, we will often refer to the image $\imm(\vertices)$ of the generic map simply as $\pos$.
For a discrete surface $\pos \in \vertexSpace$, we denote by $\len(\pos) = (\len_\edge(\pos))_{\edge\in\edges}$ its vector of edge lengths and by $ \dih(\pos) = (\dih_\edge(\pos))_{\edge\in\edges}$ its vector of dihedral angles. 
Wang \etal studied necessary and sufficient conditions that an arbitrary tuple $(\len, \dih) \in \R^{2\numE}$  is induced by a discrete surface. 
The first necessary condition is the triangle inequality, \ie
\begin{equation}
\label{eq:triangInequ}
\tag{T}
\trineq_\face(\len) > 0 \quad \text{ for all } \face \in \faces\, ,
\end{equation}
where $\trineq_f(\len) = \begin{pmatrix} l_i + l_j - l_k & l_i - l_j + l_k & - l_i + l_j + l_k \end{pmatrix}$ for a face \(\face \in \faces\) with edge lengths  \(l_i,l_j,l_k\) and the above inequality is to be understood componentwise. 
Thus, by combining the maps for all faces and extending constantly to dihedral angles we obtain a linear map $\trineq \colon \R^{2\numE} \to \R^{3 \numF}$ and we see that \eqref{eq:triangInequ} defines an open convex polytope in \(\R^{2 \numE}\).
The next set of conditions is referred to as \emph{discrete integrability conditions}. 
They ensure that we can integrate the local change of geometry induced by the lengths and angles to reconstruct the immersed discrete surface. 
In other words, the immersion is invariant with respect to the start and order of the reconstruction.
For a face \(\face \in \faces\) with (immersed) edges \(E_1, E_2, \, E_3 \in \R^3\), one defines the \emph{standard discrete frame} $\frame_\face$ as the orthogonal matrix with rows $\frac{E_1}{\lVert E_1 \rVert}$, $\frac{E_1 \times N_f}{\lVert E_1 \times N_f \rVert}$, and $N_f$, where $N_\face \in S^2$ is the unit face normal.
The normal component requires our discrete surfaces to be globally orientable which we will assume in the following. 
Then the transition between frames \(\frame_j\) and  \(\frame_i\) of adjacent faces \(\face_i, \face_j \in \faces\) and a common edge \(\edge \in \edges\) can be described by a rotation matrix \(\intRot_{ij}\) with \(\frame_j = \frame_i \intRot_{ij}\). 
This rotation decomposes into three elementary rotations, \ie
\begin{equation}\label{eq:transRotation}
\intRot_{ij} = R_2(\gamma_{e,i})R_0(-\theta_e)R_2(\gamma_{j,e}),
\end{equation}
where $R_k(\varphi) \in \SO(3)$ denotes a rotation around the $k$th standard basis vector in $\R^3$ by $\varphi \in [0,2\pi]$ and \(\gamma_{e,j}\) and \(\gamma_{i,e}\) denote the angles between the common edge and the first vector of \(\frame_i\) resp.\ \(\frame_j\).
In particular, the transition rotations are completely determined by the lengths and angles using the law of cosines. 
Now let $\vertices_0 \subset \vertices$ be the index set of interior vertices. 
Then for each \(\vertex\in\vertices_0\), which is the center of a \(\nv\)-loop of faces \(\face_0,\ldots,\face_{\nv\!-\!1}\) and edges $\edge_0, \ldots, \edge_{\nv\!-\!1}$ connected to \(\vertex\), we obtain a \emph{closing condition}. 
To this end, one chooses the frames \(\frame_0,\ldots,\frame_{\nv\!-\!1}\) such that $e_i$ always coincides with the first basis vector in $\frame_i$.
Consequently, the corresponding transition rotations simplify to \(\intRot_{ij} = R_0(\theta_i) R_2(\gamma_{j})\), where \(\gamma_{j}\) is the interior angle at \(\vertex\) in \(\face_j\) with $j = i+1$ modulo $\nv$ and $\theta_i$ is the dihedral angle at $\edge_i$.
Applying the transition property $\frame_j =  \frame_i \intRot_{ij}$ sequentially along the loop, the identity \(\frame_0 = \frame_0 \prod_{i=0}^{\nv-1} \intRot_{i,(i+1) \bmod \nv}\) must hold for immersed discrete surfaces. 
This can be phrased as the integrability condition 
\begin{equation}\tag{I}\label{eq:discreteIntegrMap}
\dint_\vertex(\len,\theta) \coloneqq \prod\limits_{i=0}^{\nv-1} R_{i,(i+1) \bmod \nv}  \overset{!}{=} \Id
\end{equation}
for the \(\nv\)-loop of faces around all \emph{interior} vertices $\vertex\in\vertices_0$. 
Note that the transition rotations in \eqref{eq:discreteIntegrMap} in fact depend on $(\len,\dih)$.

\citet{WaLiTo12} proved that the necessary conditions \eqref{eq:triangInequ} and \eqref{eq:discreteIntegrMap} are indeed sufficient (for simply connected surfaces). 
In detail, their discrete fundamental theorem of surfaces reads:
\emph{If \((\len,\dih) \in \R^{2\numE}\) satisfies \eqref{eq:triangInequ} and \eqref{eq:discreteIntegrMap}, there exists $\pos \in \vertexSpace$ (unique up to rigid body motions) such that $\len(\pos) = \len$ and $\dih(\pos) = \dih$.} 
They also extended this to non-simply connected surfaces but to simplify the exposition, we restrict ourselves to the simply connected case.

\section{The NRIC manifold of edge lengths and dihedral angles}\label{sec:nric}

In this section, we consider the Nonlinear and Rotation-Invariant Coordinates (NRIC) given as a vector $\z=(l_\edge,\theta_\edge)_{\edge \in \edges}\in\R^{2\numE}$ that lists all the edge lengths and dihedral angles of a discrete surface. Using the integrability conditions, we describe the manifold of discrete surfaces as a submanifold of $\R^{2\numE}$ and derive a scheme for computing its tangent spaces.

\paragraph*{NRIC manifold}
We consider the map 
\begin{equation}\label{eq:projection}
\projZ \colon  \vertexSpace \to \R^{\numE} \times \R^{\numE}\, , \quad  \pos \mapsto  \left(  \len(\pos),  \dih(\pos) \right) 
\end{equation}
that associates to any discrete surface the vector stacking its edge length and dihedral angles
The image of \eqref{eq:projection} describes the submanifold 
\begin{equation}\label{eq:nricDefinition}
\manifold := \projZ(\vertexSpace) =  \{\z \in \R^{2\numE} \mid \exists \pos \in \vertexSpace \colon \projZ(\pos) = \z \}.
\end{equation}
of $\R^{2\numE}$ that we call the NRIC manifold.

\paragraph*{Implicit description}
In the following, we will use the conditions \eqref{eq:discreteIntegrMap} and \eqref{eq:triangInequ} to derive an implicit description of \(\manifold\).
Directly using condition \eqref{eq:discreteIntegrMap} leads to nine scalar constraints per vertex, which is a redundant description since $\SO(3)$ is a three-dimensional manifold.
Instead, we will introduce a reformulation using unit quaternions as an equivalent representation of spatial rotations.

To this end, let us first briefly recall the necessary basics of quaternions and their relation to spatial rotation such that this section is self-contained, for a detailed treatment we refer to standard textbooks such as \cite{Ha06}.
Quaternions can be understood as an extension of the complex numbers and are generally represented as \(q = a+b \boldsymbol{i} + c \boldsymbol{j} + d \boldsymbol{k} \), where \(a,b,c,d \in \R\) and \(\boldsymbol{i}, \boldsymbol{j}, \boldsymbol{k}\) are the so-called \emph{quaternion units}.
These units fulfill the fundamental identity \(\boldsymbol{i}^2 = \boldsymbol{j}^2 = \boldsymbol{k}^2 = \boldsymbol{ijk} = -1 \), from which the general multiplication of quaternions can be defined via distributive and associative law and thus quaternions form a noncommutative division ring \(\mathbb{H}\).
In this context, $a$ is called the real part of $q$ and $b,c,$ and $d$ the vector part, for which we also write $\mathrm{vec}(q) = (b,c,d) \in \R^3$.
Unit quaternions are those for which the product with their conjugate \(\bar q \coloneqq a-b \boldsymbol{i} - c \boldsymbol{j} -  d \boldsymbol{k} \) is one, \ie $q \bar q = a^2 + b^2 + c^2 + d^2 = 1$.
Points in three-dimensional space \(p \in \R^3\) can be identified with quaternions having vanishing real part, \ie we write \(p = p_1 \boldsymbol{i} + p_2 \boldsymbol{j} + p_3 \boldsymbol{k}\).
Now, given a rotation \(Q\) around the unit vector \(u \in \R^3\) by angle \(\varphi \in [0, 2 \pi)\) we can define a corresponding \emph{unit} quaternion 
\[q(u,\varphi) \coloneqq \cos \frac{\varphi}{2} + \left(u_1 \boldsymbol{i} + u_2 \boldsymbol{j} + u_3 \boldsymbol{k}\right) \sin \frac{\varphi}{2}.\]
Then one can verify that for any \(p \in \R^3\), the conjugation $qpq^{-1}$ with $q$ results in the rotated point \(Qp\). 
The quaternion $-q(u, \varphi)$ would lead to the same rotation, thus the quaternions form a double covering of \(\SO(3)\).  
Furthermore, investigating this conjugation one realizes that the composition of two rotations given as unit quaternions \(q_1, q_2 \in \mathbb{H}\) is given by their product \(q_1q_2\) and hence this correspondence is a homomorphism between \(\SO(3)\) and the unit quaternions.

Turning to the reformulation of the integrability conditions \eqref{eq:discreteIntegrMap}, recall that we needed rotations around the $0$th and $2$nd basis vector in $\R^3$ for which we now introduce the corresponding quaternions
\begin{equation}
	q_0(\varphi) \coloneqq \cos \frac{\varphi}{2} + \boldsymbol{i} \sin \frac{\varphi}{2}, \quad q_2(\varphi) \coloneqq \cos \frac{\varphi}{2} + \boldsymbol{k} \sin \frac{\varphi}{2} \quad \text{ for } \varphi \in [0, 2 \pi).
\end{equation}
Then we identify the simplified transition rotation \(\intRot_{ij} = R_0(\theta_i) R_2(\gamma_{j})\) from before with the quaternion 
\begin{equation}
	q_{ij} \coloneqq q_0(\theta_i)\, q_2(\gamma_{j}),
\end{equation}
where again \(\face_0,\ldots,\face_{\nv\!-\!1}\) and $\edge_0, \ldots, \edge_{\nv\!-\!1}$ are the \(\nv\)-loops of faces resp.\ edges connected to $\vertex$, \(\gamma_{j}\) is the interior angle at \(\vertex\) in \(\face_j\) with $j = i+1$ modulo $\nv$, and $\theta_i$ is the dihedral angle at $\edge_i$. 
To finally reformulate the condition \eqref{eq:discreteIntegrMap}, we need to deal with the ambiguity introduced by the double covering, \ie that the identity rotation is represented by \(q = \pm 1\). 
However, we see that in both cases the vector part \(\mathrm{vec}(q) \in \R^3\) is zero, which is indeed for unit quaternions already a sufficient condition to be plus or minus one.
Then we use this alternative characterization of the identity rotation to formulate the \emph{quaternion integrability conditions} as
\begin{equation}\tag{I\textsubscript{q}}\label{eq:quatIntegrMap}
\qint_\vertex(\len,\theta) \coloneqq \mathrm{vec}\left(\prod\limits_{i=0}^{\nv-1} q_{i,(i+1) \bmod \nv}\right)  \overset{!}{=} 0
\end{equation}
for the \(\nv\)-loop of faces around all interior vertices $\vertex\in\vertices_0$. 

Now, we can rewrite the manifold defined in \eqref{eq:nricDefinition} as 
\begin{align}\label{eq:nricDefinitionImplicit}
\manifold = \big\{\z \in \R^{2\numE} \,\big|\, \trineq(z) > 0,\, \qint(z) =0 \big\}\, .
\end{align} 
Here we have collected all constraints in a vector-valued functional $\qint \colon \R^{2\numE} \to \R^{3|\vertices_0|}$ with $\qint = (\qint_\vertex)_{\vertex \in \vertices_0}$. 
Obviously, $\qint_\vertex$ depends solely on the edge lengths of the adjacent faces of $\vertex$ and the dihedral angles at edges centered at  $\vertex$.
Given $\z \in \R^{2\numE}$ with $\qint(\z) = 0$ one can easily reconstruct vertex coordinates $\pos \in \vertexSpace$ with $\projZ(\pos) =\z$. 
For a robust and stable reconstruction for $\qint(\z) \neq 0$, we refer to Section~\ref{sec:recon}.


\paragraph*{Tangent space}
The implicit formulation \eqref{eq:nricDefinitionImplicit} consists of the triangle inequalities defining an open convex polytope and of the nonlinear integrability conditions, which define a lower-dimensional, differential structure on $\manifold$. 
Therefore, we can derive an implicit description of its tangent space solely based on $\qint$.
In detail, for $\z\in\manifold$ the tangent space is given by
\begin{align} \nonumber
T_\z \manifold &= \mathrm{ker}\, D\qint(\z) \coloneqq \{ w \in \R^{2\numE} \, | \, D\qint(\z) w = 0 \}\,,
\end{align}
where $D \qint(\z)$ is a matrix in $\R^{3|\vertices_0|, 2\numE}$. 
Partial derivatives of $\qint_\vertex$ are given by the chain rule as
\begin{align}\label{eq:intFirstDeriv}
\partial_{\z_k} \qint_\vertex(\z) &= \mathrm{vec} \left(\sum_{i=0}^{\nv-1} q_{01}(\z) \ldots \partial_{\z_k} q_{i,i+1}(\z)   \ldots  q_{\nv-1,0}(\z)\right),
\end{align}
where the partial derivatives of a quaternion-valued map are to be understood componentwise as for vector-valued maps.
The gradient of $\qint_\vertex$ can be computed with $O(\nv)$ cost and is sparse.
It has only $O(\nv)$ non vanishing entries, \ie $\partial_{\theta_\edge} \qint_\vertex \equiv 0$ if $\vertex$ is not a vertex of the edge $\edge$ and $\partial_{\len_\edge} \qint_\vertex \equiv 0$ if the edge $\edge$ is not an edge of a triangle with vertex $\vertex$.
We provide details on the gradient computation as well as an implementation in terms of a \textsc{Mathematica} notebook in the supplementary material. 
\begin{figure}[ht]
	\centering
	\includegraphics[width=0.24\columnwidth]{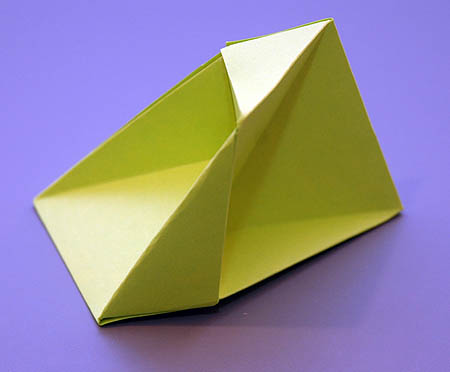}
	\includegraphics[width=0.24\columnwidth]{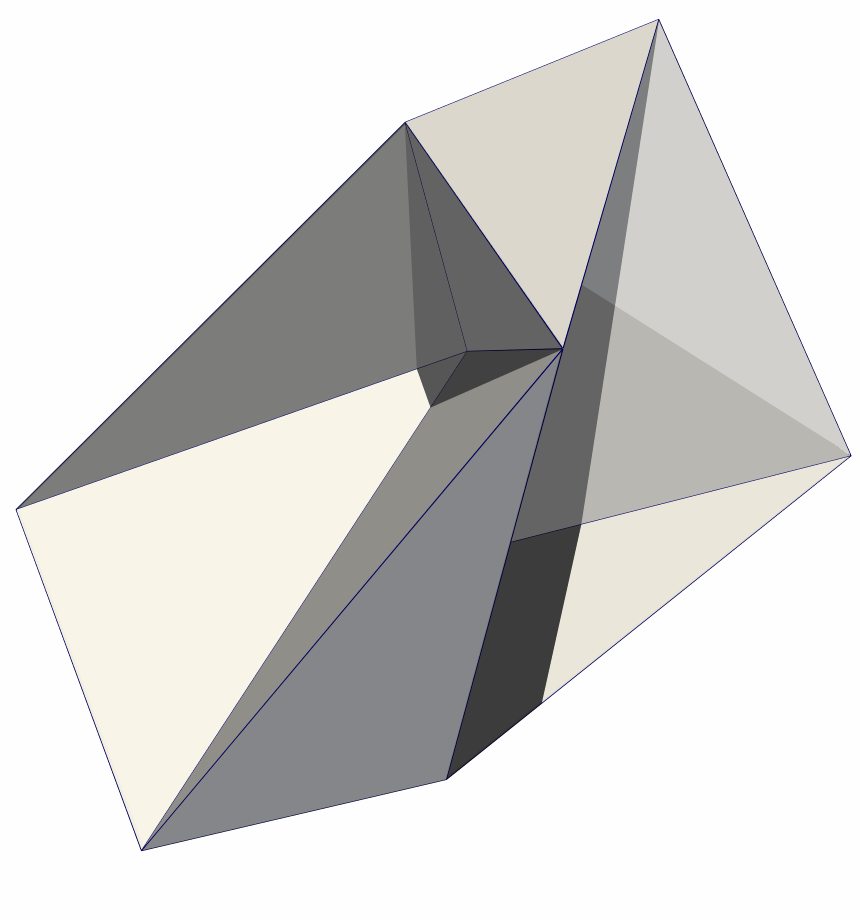}
	\includegraphics[width=0.24\columnwidth]{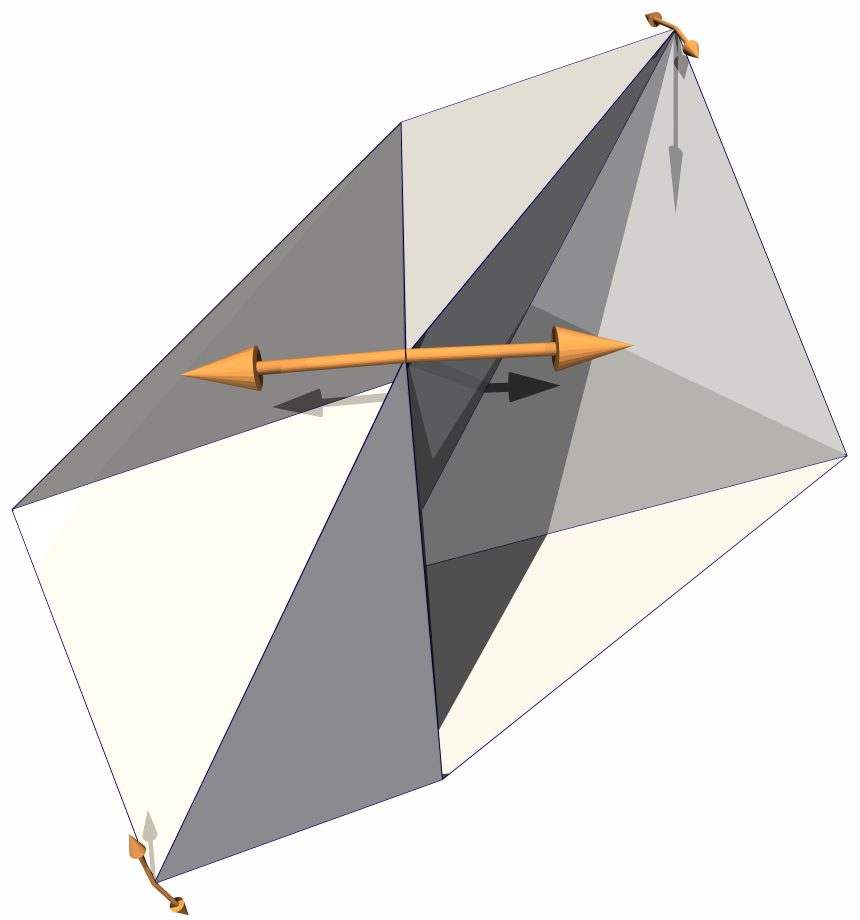}
	\includegraphics[width=0.24\columnwidth]{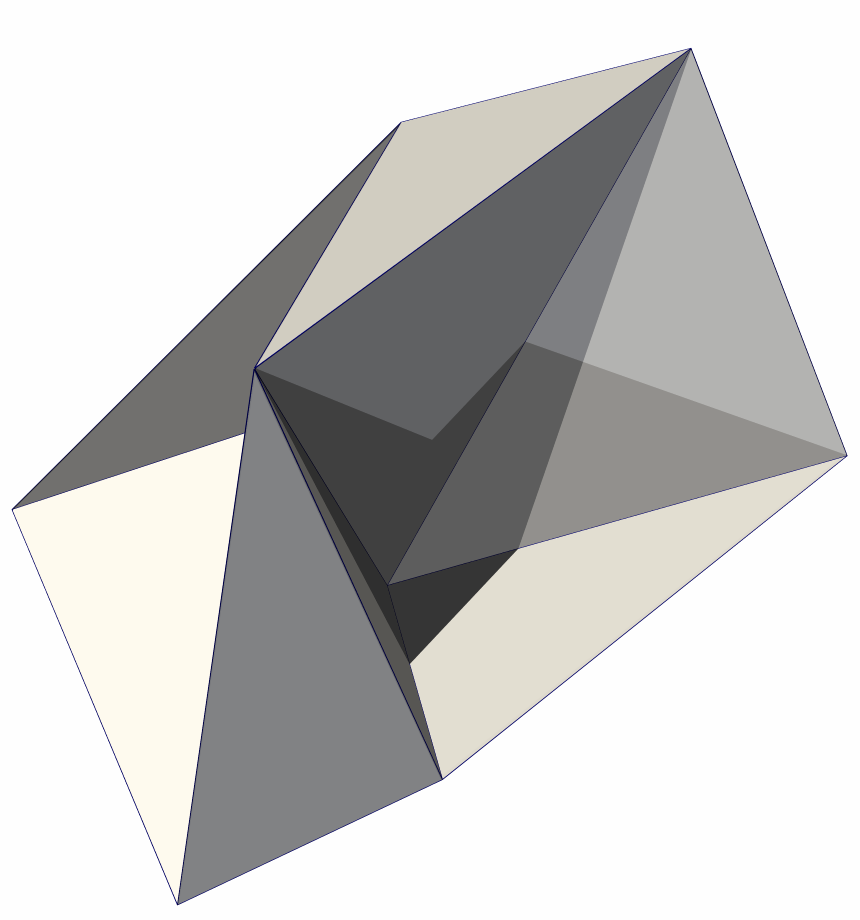}
	\caption{The tangent space reveals an infinitesimal isometric variation at the classical Steffen's polyhedron (middle). 
		Indeed, extrapolating in this positive (left) resp. negative (right) direction (solely in the $\theta$ component) allows for isometric deformations. 
		The extrapolation is implemented via an incremental addition of the infinitesimal isometric variation coupled with a back projection onto $\manifold$. 
		See also video in supplementary material for an animation. }
	\label{fig:steffens}
\end{figure}

\begin{figure}[ht]
	\centering
	\includegraphics[width=0.19\columnwidth]{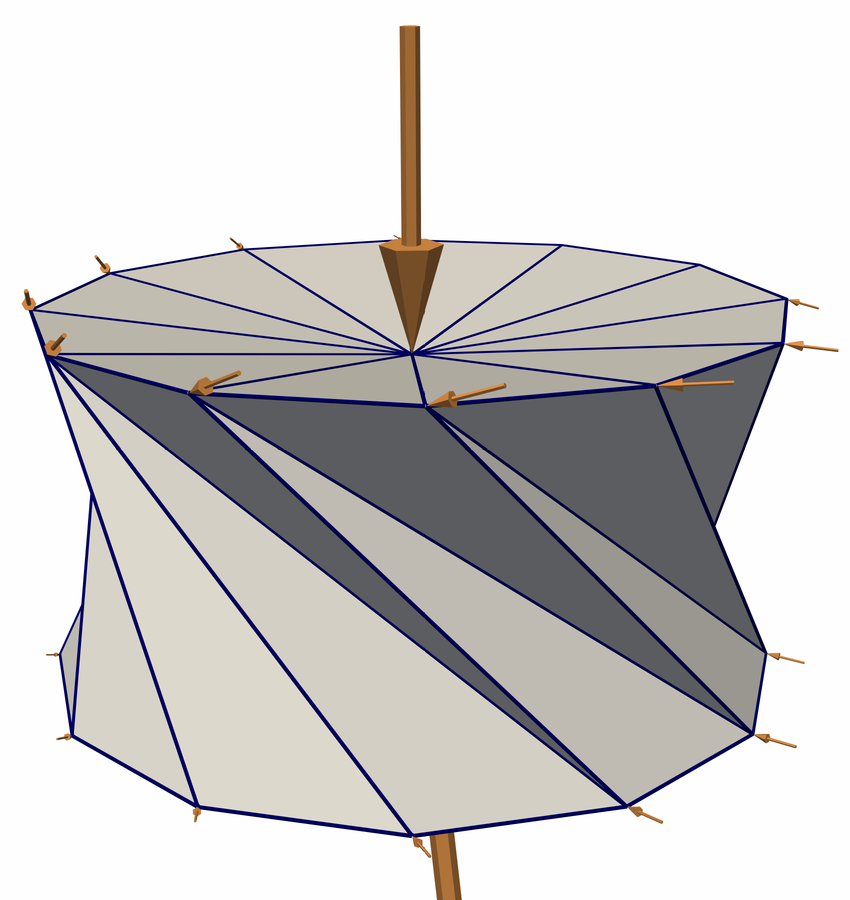}
	\includegraphics[width=0.19\columnwidth]{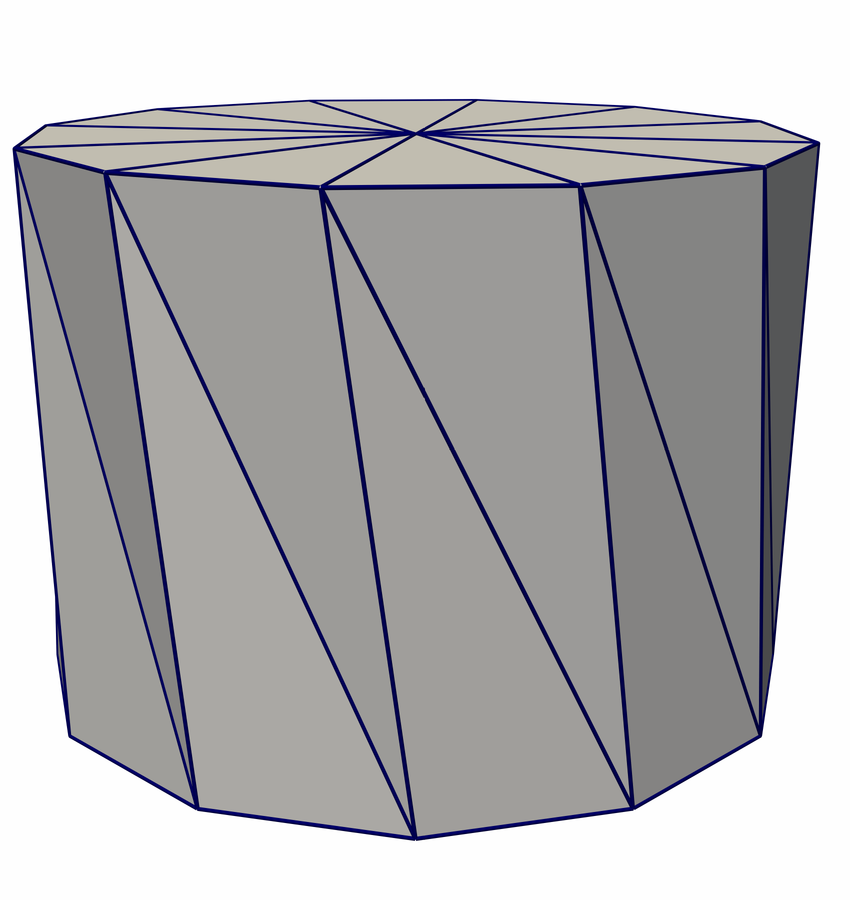}
	\includegraphics[width=0.19\columnwidth]{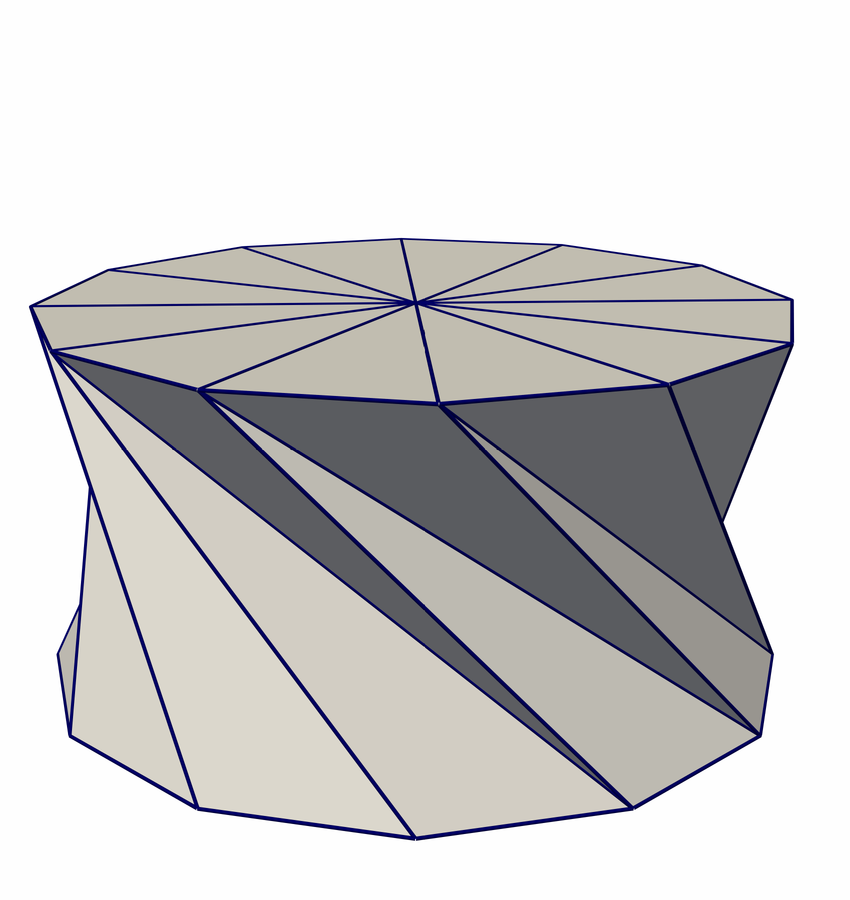}
	\includegraphics[width=0.19\columnwidth]{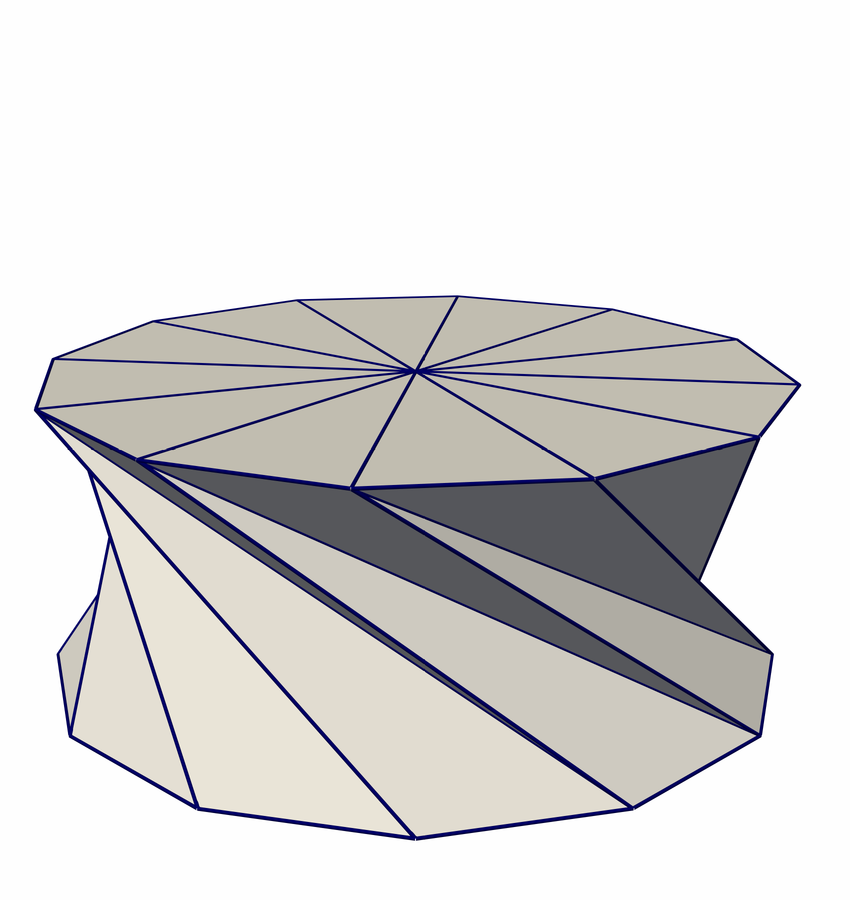}
	\includegraphics[width=0.19\columnwidth]{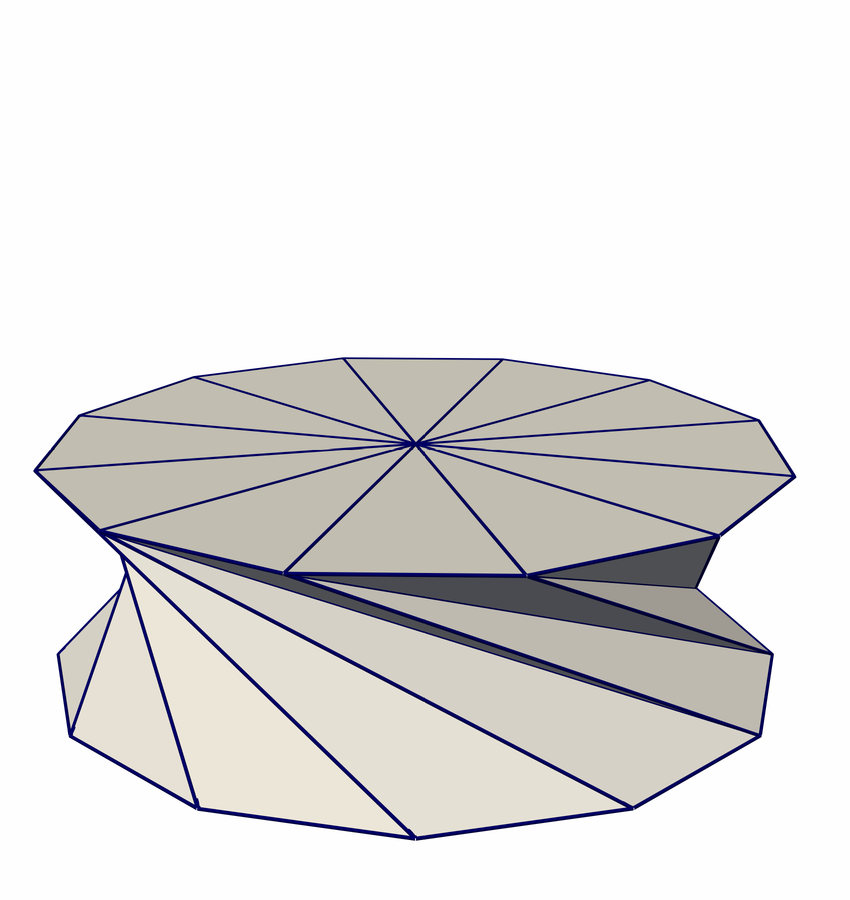}\\  
	\parbox{0.19\linewidth}{\hfill}
	\parbox{0.19\linewidth}{\centering $\dih^\ast\!= 0.75\cdot\bar \dih$}
	\parbox{0.19\linewidth}{\centering reference}
	\parbox{0.19\linewidth}{\centering $\dih^\ast\! = 1.1 \cdot \bar \dih$ }
	\parbox{0.19\linewidth}{\centering $\dih^\ast \!= 1.2 \cdot \bar \dih$}
	\caption{
		Top: Almost isometric compression of an Origami cylinder as depicted in \citet{BoVoGoWa16} with the only infinitesimal isometric variation (left). 
		Optimizing \eqref{eq:quadEnergy} on NRIC manifold with $\delta\!=\!0$ and hard constraints on target angles $\dih^\ast$ along the upper horizontal edges (relative to the reference angle $\bar\dih = 2.257$) leads to non-isometric deformations with as small as possible edge length distortion (from left to right $4\%$, $0.3\%$, and $0.4\%$ average change of edge length). 
		See also supplementary video.
	}
	\label{fig:origami}
\end{figure}

To illustrate the NRIC manifold and its tangent spaces, we will for the remainder of the section discuss an immediate application. 
With the tangent space at hand, one can verify the \emph{infinitesimal rigidity} of a discrete surface with NRIC $\z\in \manifold$.
In fact, a necessary condition for the existence of continuous one-parameter families of isometric deformations starting at $\z$ is the existence of an infinitesimal isometric variation $w\in T_\z \manifold$ with $P_\len w = 0$ and $w\neq 0$, where $P_\len$ is the projection onto the length component, \ie  $P_\len (\len,\dih) = \len$, see for example \autoref{fig:steffens}. 
Note, however, that this is surely not a sufficient condition, which we can also observe in \autoref{fig:origami}.
Thus, we simply verify if the kernel of $D \qint(\z)$ has a non-trivial intersection with the $\dih$ subspace of $\R^{2\numE}$, namely the kernel $\ker P_\len$.
This intersection is given by $\ker \left( \begin{array}{c|c} B_{T_\z \manifold} & B_{\theta} \end{array}\right)\,$, 
where $B_{T_\z \manifold}$ is a matrix whose columns form a orthonormal basis of $T_\z\manifold$ and 
$B_{\theta}$ is the canonical basis of $\ker P_\len$.
We compute a singular value decomposition (SVD) of this matrix and evaluate the smallest singular value $\lambda_0$. 
If $\lambda_0 = 0$, then there exists an infinitesimal isometric variation. 
Otherwise, the singular value provides a quantitative measure for the lack of such an infinitesimal isometric variation. 
In \autoref{fig:origami}, we show that the Origami cylinder considered by \citet{BoVoGoWa16} does not allow for an isometric deformation path which leads to a compression by folding. 
Indeed, the only nontrivial infinitesimal isometric variation is indicated by arrows (top, right).
However, there is no nontrivial family of isometric deformation with this shape as the initial shape.
As discussed by \citet{BoVoGoWa16} the experimental paper deformation (top, left) is not isometric. 
This is reflected by our criteria for infinitesimal isometric variations when we additionally enforce the dihedral angles on the upper and lower plate to remain constant which leads to $\lambda_0=0.015$ clearly indicating the nonexistence of such a variation.


\section{Nonlinear energy and geometry of the NRIC manifold}
\label{sec:energy}
So far, we have introduced a differential structure on $\manifold$.
Going forward, it will be essential to additionally consider an elastic deformation energy $\W$ between different NRIC as it provides a dissimilarity measure on $\manifold$.
Although different choices for $\W$ are possible, we will primarily focus on the hyperelastic deformation energy from \citet{HeRuSc14}. 
To this end, we reformulate this energy in NRIC to define a physically-motivated Hessian structure on $\manifold$.
This will be a straightforward undertaking which underlines our claim that NRIC are a natural choice for computing deformations.
In particular, we will see that the local injectivity constraints inbuilt in this energy allow us to replace the triangle inequalities and thus reduce the number of constraints. 
For comparison reasons, we will finally consider a simple quadratic deformation energy as it has been used in \cite{FrBo11}.

Based on models from mathematical physics, the hyperelastic energy used in \cite{HeRuSc14} consists of two separate contributions, i.e.
\begin{equation}
	\label{eq:genericEnergy}
	\mathcal{W} = \Wmem + \delta^2\, \Wbend.
\end{equation}
From a physical point of view, the first term $\Wmem$ will measure the stretching of edges and triangles, i.e. local \emph{membrane distortions}. 
Likewise, the second term $\Wbend$ will measure the difference in bending between triangles, i.e. local \emph{bending distortions}.
In particular, the global weight $\delta$ represents the thickness of a thin elastic material represented by the discrete surface.

In the following, we will investigate separately how the membrane and bending energy introduced in \cite{HeRuSc14} can be reformulated in NRIC. 
Note that the membrane energy has originally been proposed in \cite{HeRuWa12} whereas the bending energy has been taken from \cite{GrHiDeSc03}. 

\paragraph{Membrane energy}
Let \(\imm\in\vertexSpace\) be a discrete surface and \(\immEx\colon \underlyingSpace \to \R^3\) the corresponding continuous, piecewise linear map that interpolates the vertices. 
The derivative of \(\immEx\) is constant over each triangle, and, since \(\imm\) is a generic map, the derivative has full rank. 
This implies that \(\immEx\) induces a metric $\dfFF$ (also called first fundamental form) on $\underlyingSpace$. 
This metric is defined in the interior of the faces and along the edges. It enables measuring the length of arbitrary curves in $\underlyingSpace$ and makes $\underlyingSpace$ a metric space.  
Two discrete surfaces, \(\imm\) and \(\immAlt\), induce two different metrics, $\dfFF$ and $\widetilde \dfFF$, on $\underlyingSpace$. 
The metric distortion tensor $\dDisTen[\imm , \immAlt]$ is defined as the symmetric tensor that at any point in the interior of a triangle satisfies \[\dfFF(\dDisTen[\imm ,\immAlt]v,w) = \widetilde\dfFF(v,w)\] for any pair $v,w$ of tangential vectors. 
The membrane energy evaluates the trace and the determinant of $\dDisTen[\imm, \immAlt]$. 

For our purpose, it is essential to be able to evaluate the distortion tensor for discrete surfaces 
\begin{wrapfigure}{r}{0.2\columnwidth}
    \vspace{-1.3em}
    \hspace{-2em}
    \small 
    \def\svgwidth{100pt}
    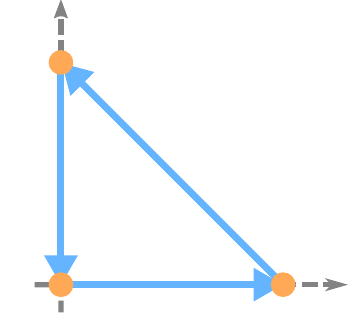
\end{wrapfigure}
given by their NRIC $\z$ and $\tilde \z$ directly without having to reconstruct vertex positions first. 
In the following, we derive an explicit formula for $\dDisTen[\z ,\tilde \z]$. 
For discrete surfaces, the metric and the distortion tensor are constant for every triangle. 
Consider an arbitrary triangle \(\face\) in $\R^3$. We parametrize \(\face\) with an affine map $\phi\colon t \to  \face$, where $t$ is the right angled triangle in $\R^2$ shown in the inset figure. 
The standard basis $b_1,b_2$ of $\R^2$ agrees with second edge and the negative of the first edge of $t$. 
Then $\d\phi(b_1)=E_2(f)$ and $\d\phi(b_2)=-E_1(f)$, where $E_1(f),E_2(f)$ denote the edge vectors of $\face$. 
Thus, the metric on $t$ induced by $\phi$ is 
\begin{equation} \label{eq:discFirstFundForm}
	\dfFF\vert_\face =
\begin{pmatrix}
	\lVert \d\phi(b_1) \rVert^2 &  \langle \d\phi(b_2), \d\phi(b_1)\rangle \\
	 \langle \d\phi(b_2), \d\phi(b_1)\rangle & \lVert \d\phi(b_2) \rVert^2
	\end{pmatrix}
=	
	\begin{pmatrix}
	\lVert E_2(f) \rVert^2 & - \langle E_1(f), E_2(f)\rangle \\
	- \langle E_1(f), E_2(f)\rangle & \lVert E_1(f) \rVert^2
	\end{pmatrix},
\end{equation}
The entries of the metric can be expressed in terms of the length of the edges of $\face$.
The diagonal entries are the squared length of the second and the first edge. 
The off-diagonal entries are given by scalar products of edge vectors and from linear algebra we recall that for two vectors \(v, w \in \R^3\) we have \(\langle v, w\rangle = \lVert v \rVert \, \lVert w \rVert \cos(\gamma) \), where \(\gamma\) is the angle between \(v\) and \(w\).
In our case, this is the interior angle of a triangle which can be computed from its edge lengths by the law of cosines.
For two NRIC $\z$ and $\tilde \z$, we can use the formula to compute the metrics $\dfFF\vert_\face$ and $\widetilde\dfFF\vert_\face$ for every $\face$. 
Then, the distortion tensor is given as \(\dDisTen[\pos, \tilde \pos]\vert_\face := (\dfFF\vert_\face)^{-1}\widetilde \dfFF\vert_\face\). 
We want to note that the resulting distortion tensor depends on the chosen domain and parametrization. 
However, we consider isotropic materials for which the membrane energy depends only on the trace and determinant of the distortion tensor. 
Since the determinant and the trace are invariant under coordinate transformations, we obtain the same results independently of the chosen domain and parametrization. 
Similarly the roles of the edges could be exchanged, for example, one could consider the second and third edge. 
This would alter the parametrization and therefore yield a different distortion tensor. 
Still, the relevant quantities, the determinant and the trace of $\dDisTen$, would be the same.  
  
Having established that the distortion tensor is completely given by the NRIC of discrete surfaces, we can now adapt the membrane energy from \cite{HeRuWa12} applying a nonlinear energy density to it, which has a global minimum at the identity.
\begin{definition}[Membrane energy]
	\label{def:ltheta_membrane_energy}
	For a simplicial surface \(\simplicial\), we define the membrane energy on NRIC \(\z, \tilde \z \in \R^{2\numE}\) as
	\begin{equation}
		\label{eq:ltheta_membrane_energy}
		\W_\mem[\z, \tilde \z ] = \sum_{\face \in \faces} a_\face \cdot W_\mem(\mathcal G[\z, \tilde \z]|_\face),
	\end{equation}
	where 
	\begin{equation*}
		W_{\mem}(A) := \frac{\mu}{2}\tr A + \frac{\lambda}{4}\det A -\left(\mu+\frac{\lambda}{2}\right)\log \det A - \mu - \frac{\lambda}{4},
	\end{equation*}
	for positive material constants \(\mu\) and \(\lambda\) and $\area_\face$ is the area of $\face$ computed from edge lengths by Heron's formula. 
\end{definition}
For more explicit formulas of the energy in terms of edge lengths we refer to the appendix and for the energy's derivatives to the supplementary material.

\paragraph{Bending energy}
Next, we adapt the \emph{Discrete Shells} bending energy \cite{GrHiDeSc03} also used in \cite{HeRuSc14}.
One directly sees that expressing this energy in lengths and angles requires no further calculations, and as before we replace its primary variables by NRIC.
\begin{definition}[\emph{Discrete Shells} bending energy]
	\label{def:ltheta_bending_energy}
	For a simplicial surface \(\simplicial\), we define the Discrete Shells bending energy on NRIC \(\z, \tilde \z \in \R^{2\numE}\) as
	\begin{equation}
	\label{eq:ltheta_bending_energy}
	\W_\bend[\z, \tilde \z] = \sum_{\edge\in\edges} \frac{(\theta_e - \tilde{\theta}_e)^2}{d_\edge}l_\edge^2,
	\end{equation}
	where \(d_\edge = \frac{1}{3}(\area_\face + \area_{\face'})\) for the two faces \(\face\) and \(\face'\) adjacent to \(\edge \in \edges\), as before computed by Heron's formula.
\end{definition}

Finally, we combine the membrane and bending energy in a weighted sum.
\begin{definition}[Nonlinear deformation energy]
	\label{def:ltheta_energy}
	Let \(\simplicial\) be a simplicial surface and let \(\z, \tilde \z \in \R^{2\numE}\) be two NRIC.
	The nonlinear deformation energy is defined by 
	\begin{equation} \label{eq:nonlinearEnergy}
	\W_{nl}[\z, \tilde \z] = \W_\mem[\z, \tilde \z ] + \delta^2\, \W_\bend [\z, \tilde \z ],
	\end{equation}
	where $\W_\mem$ is the membrane energy from \autoref{def:ltheta_membrane_energy}, $\W_\bend$ is the bending energy from \autoref{def:ltheta_bending_energy}, and \(\delta\) represents the thickness of the material.
\end{definition} 

\paragraph{Relationship with triangle inequalities}
One essential property of the membrane energy is that it allows us to control local injectivity via the built-in penalization of volume shrinkage, \ie we have \(W_\mem(\mathcal G[\z, \tilde \z]|_\face) \to \infty\) for \(\tilde{\area}_\face \to 0\).
To see this, we recognize that $\det \dDisTen[\z, \tilde \z]|_\face = (\det \dfFF\vert_\face)^{-1} \det \widetilde \dfFF \vert_\face = \area_\face^{-2}\,\tilde \area_\face^2$ and hence $-\log \det \dDisTen[\z, \tilde \z]|_\face \to \infty$ when $\tilde{\area}_\face$ goes to zero.
This control over the local injectivity also has consequences for the consideration of the triangle inequalities.
Because of it, we also have that the energy diverges, \ie \(W_\mem(\mathcal G[\z, \tilde \z]|_\face) \to \infty\) if one of the components of \(\trineq_\face(l)\) approaches zero meaning that we get close to violating one of the triangle inequalities.
Especially, we set \(W_\mem(\mathcal G[\z, \tilde \z]|_\face) = \infty\) if \(\trineq_\face(l) > 0\) does not hold.
This allows us to characterize the NRIC manifold \(\manifold\) by 
\begin{equation}
\label{eq:manifoldViaEnergy}
\manifold = \big\{\z \in \R^{2\numE} \,\big|\, \W_{nl}[\z^\ast, z] < \infty \text{ for a fixed } \z^\ast \in \manifold ,\, \qint(z) =0 \big\}\, ,
\end{equation}
avoiding the explicit dependence on the triangle inequalities \autoref{eq:triangInequ} we had before.
Note, however, that the integrability conditions \autoref{eq:quatIntegrMap} are still necessary as finite energy does not guarantee their attainment.
The characterization \autoref{eq:manifoldViaEnergy} will be helpful later on to devise efficient numerical schemes for solving variational problems on $\manifold$.

\paragraph{Quadratic model} 
Previously, \citet{FrBo11} used a quadratic deformation model for NRIC, \ie they considered the weighted quadratic energy
\begin{equation} \label{eq:quadEnergy}
\W_q[\z,\z^\ast] = \edgeSum \memw_\edge \lVert \len_\edge -  \len_\edge^\ast \rVert^2 + \thickSqr \,  \edgeSum \benw_\edge \lVert \dih_\edge - \dih_\edge^\ast \rVert^2\, .
\end{equation}
In fact, almost the same model has been used in \cite{GrHiDeSc03} to define the \emph{Discrete Shells} energy for physical simulations based on nodal positions.
The weights $\memw = (\memw_\edge)_\edge$ and $\benw = (\benw_\edge)_\edge$ can be chosen in different ways. 
Typically, they are computed from edge lengths $\len_\edge = \len_\edge(\bar z)$ 
and areas $d_\edge = d_\edge(\bar \z)$ associated with edges and defined on some representative reference configuration $\bar\z\in\R^{2\numE}$. 
For example, the authors in \cite{GrHiDeSc03, FrBo11} set in a related context $\memw_\edge =  \len_\edge^{-2}$ and $\benw_\edge =  \len_\edge^2 \,  d_\edge^{-1}$, whereas \citet{HeRuSc16} have chosen $\memw_\edge =  d_\edge^{\vphantom{-2}}  \len_\edge^{-2}$, for $\edge \in \edges$.
Here the (physical) parameter $\thickSqr > 0$ trades the impact on length variations off against angle variations and can be considered as the squared thickness of the material as before.
This quadratic energy has no inbuilt control over the local injectivity of the deformation and hence does not allow a characterization without explicit dependence on the triangle inequalities as in \autoref{eq:manifoldViaEnergy}.
We found that in many of our examples this decreased the numerical accessibility and increased the needed number of iterations and runtimes.
Nevertheless, as demonstrated by \citet{FrBo11}, it often leads to natural-looking deformations and we will consider it in some of our examples.

\paragraph{Riemannian metric}
For each $\z\in \manifold$, a Riemannian metric $\metric_\z$ is a symmetric, positive definite quadratic form on the tangent space $T_\z \manifold$ measuring the cost of an infinitesimal variation in tangential direction. 
In our context a tangential vector $w\in T_\z \manifold$ splits into two components $w= (w_\len, w_\dih)$, where $w_\len$ is the variation of edge lengths and $w_\dih$ the variation of dihedral angles. 
Following Rayleigh's paradigm, $\frac12$ times the Hessian of an elastic deformation energy can be considered as a Riemannian metric on the space of discrete surfaces if it is positive definite. 
Precisely, we obtain the metric for tangent vectors $v,w \in \R^{2\numE}$ via
\begin{align}\label{eq:metric}
\metric_{\z}(v,w) = v^T\left(\frac12 \Hess \W[\z,\z] \right)w\, .
\end{align}
As investigated in \cite{HeRuSc14}, this is true for the energy defined in \autoref{eq:nonlinearEnergy} with the choice of membrane and bending energies made above.
Furthermore, it holds for the quadratic energy \autoref{eq:quadEnergy} if we choose all weights to be positive.
With the metric at hand, one can define the Riemannian distance on $\manifold$ and compute for instance shortest geodesic curves, \cf Section~\ref{sec:results}.

\section{Variational problems on the manifold}\label{sec:varProblems}
The quest for geometrically optimal, discrete surfaces often leads to variational problems. 
However, in many applications, the corresponding objective functional can naturally be formulated in our coordinates, thus on the NRIC manifold \eqref{eq:nricDefinitionImplicit}, and its first and second variation can be computed easily.
To this end, one aims at solving a constrained optimization problem, \ie given an objective functional $\energy \colon \R^{2\numE} \to \R$ the task is to
\begin{equation}
	\begin{aligned}
	& \underset{\z\, \in\, \R^{2\numE}}{\text{minimize}} 
	& & \energy(\z) \\
	& \text{subject to}
	&  & \qint_\vertex(\z) = 0 \text{ for each } \vertex \in {\vertices_0}, \\
	& 
	&  & \trineq_\face(\z) > 0 \text{ for each } \face \in {\faces}.
	\end{aligned}
	\tag{OPT}
	\label{eq:genericOpt}
\end{equation}
Due to non-convexity of the objective, in general, there is no guarantee for a unique, global minimizer for the optimization problem.

\begin{figure}[ht]
	\begin{minipage}[c]{0.56\columnwidth}
		\includegraphics[trim=70 25 75 25, clip,width=0.48\columnwidth]{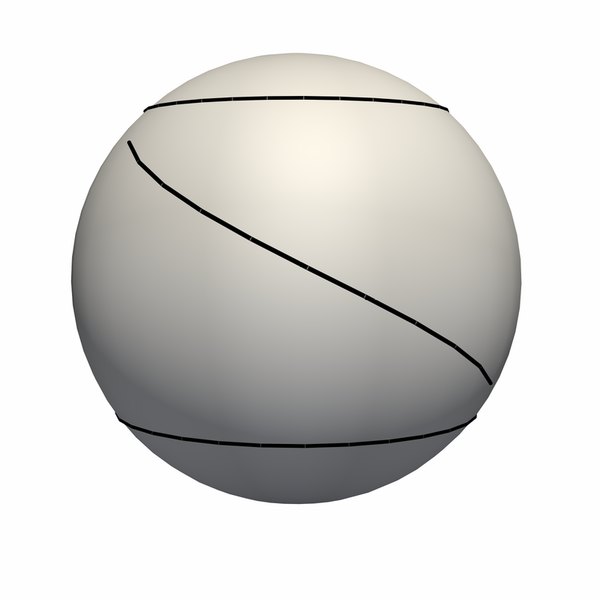}
		\includegraphics[trim=75 25 75 25, clip,width=0.48\columnwidth]{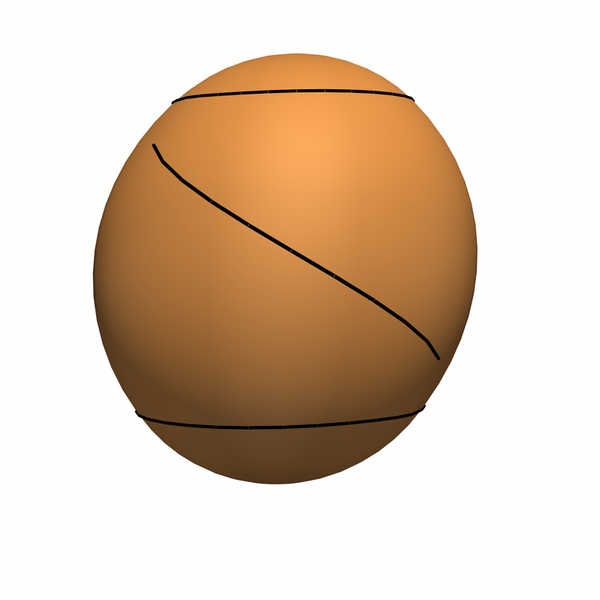}\\
		\includegraphics[trim=75 0 75 50, clip,width=0.48\columnwidth]{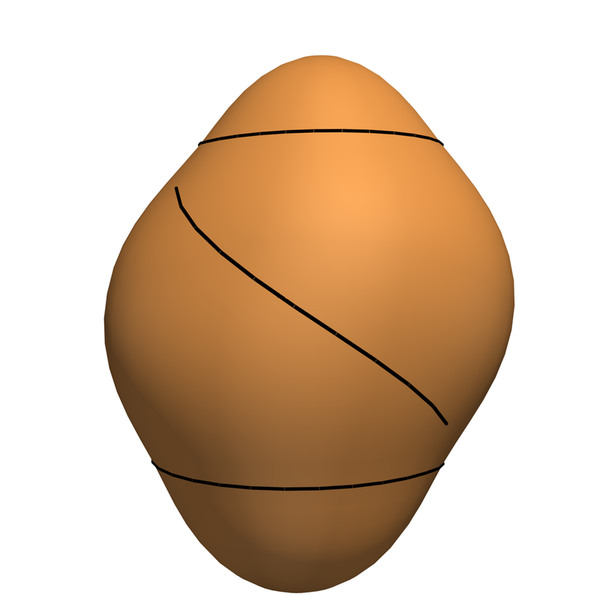}
		\includegraphics[trim=75 25 75 25, clip,width=0.48\columnwidth]{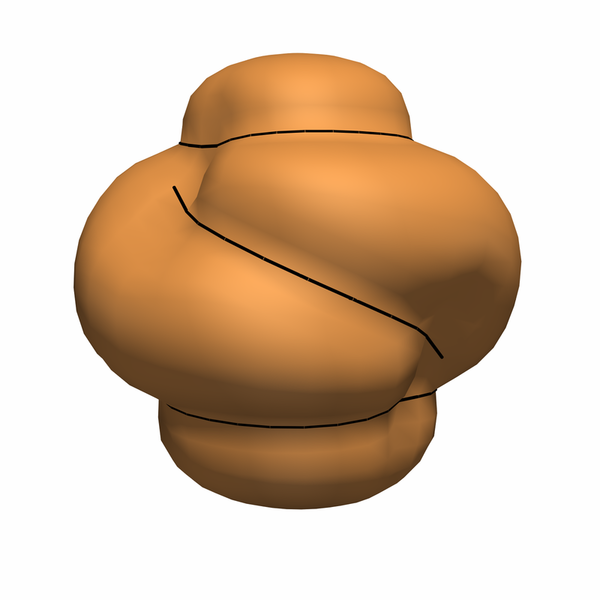}		
	\end{minipage}
	\begin{minipage}[c]{0.39\columnwidth}
		\includegraphics[trim=100 0 80 0, clip,width=0.48\columnwidth]{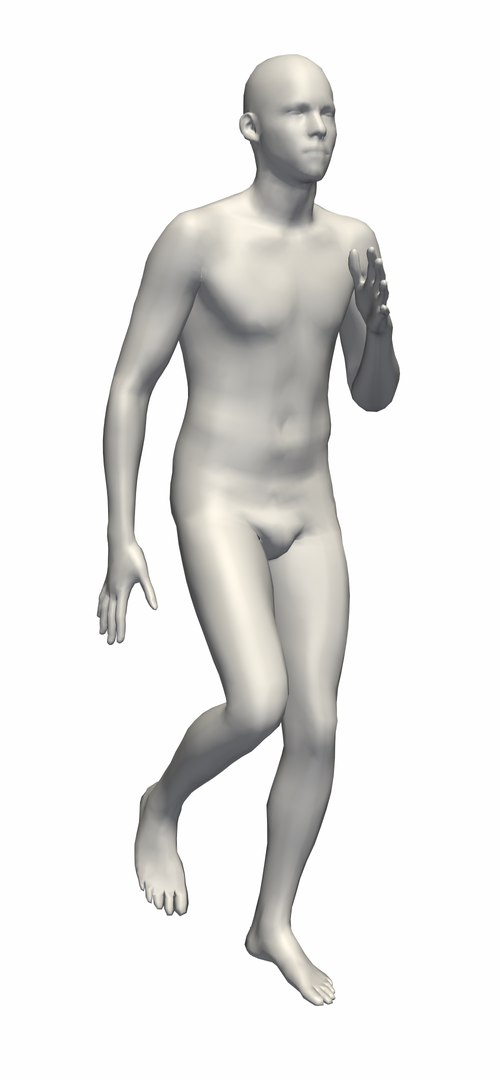}
		\includegraphics[trim=100 0 80 0, clip,width=0.48\columnwidth]{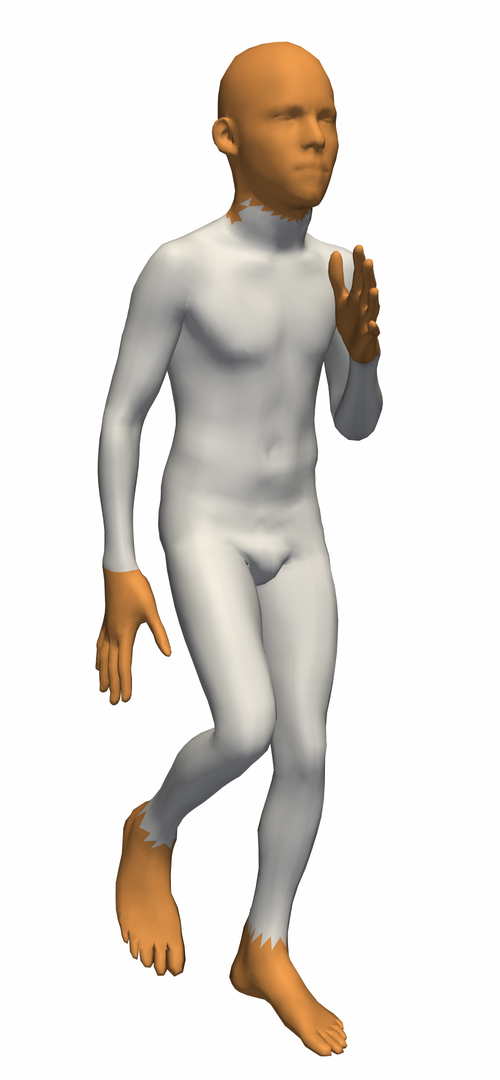}
	\end{minipage}
	\caption{Left: Unit sphere (grey) with black constraint curves to be shortened by means of equality constraints on edge lengths along with different results for varying bending parameter $\thickSqr= 10^{-\{0,2,4\}}$ in the nonlinear  objective \autoref{eq:nonlinearEnergy}. 
		Right: Same experiment but with (orange) constraint areas where the target edge lengths were increased by $30\%$ to simulate a brushing tool.}
	\label{fig:strangulation}
\end{figure}

A simple example of an objective functional $\energy$ is given by the dissimilarity to some given $\z^\ast$ on the linear space $\R^{2\numE}$ measured by the deformation energy, \ie $\energy(\z) = \W[\z^\ast,z]$, where $\W$ is an elastic deformation energy as discussed in the previous section. 
For example, in \autoref{fig:strangulation} we have used the nonlinear energy defined in \autoref{eq:nonlinearEnergy} along with coordinate constraints on a certain subset of edge lengths to simulate a ``constriction'' of a sphere along curves or creating cartoon-like characters by inflating for instance hands and feet (\cf \cite{KrNiPo14}). 

\paragraph{Ensuring triangle inequalities}
One crucial problem we encounter when we try to solve \autoref{eq:genericOpt} are the triangle inequalities which lead to $3\numF$ inequality constraints causing the problem to be computationally expensive.
Therefore, we aim for an approach to deal with them efficiently rooted in our geometric setup from \autoref{sec:nric} and \autoref{sec:energy}.
We achieve this by a modified line search.
First, recall that the set \(\trineqSet = \left\{ \z \in \R^{2\numE} \mid \trineq(\z) > 0 \right\} \) defines an open connected subset of $\R^{2\numE}$.
Therefore, if we start with an initial point $\z^0$ fulfilling the triangle inequalities we only have to ensure that every iterate remains in the set.
Hence, in a line search method where we search for a new iterate $\z^{k+1}$ along a direction $d^k$ we have to restrict this search to $\trineqSet$.
We accomplish this using backtracking, \ie reducing the stepsize $\beta^k$ until $\z^{k+1} = \z^k + \beta^k d^k \in \trineqSet$ holds.
In implementations, this can easily achieved by setting $\qint_\vertex(\z) = \infty$ if $\trineq_\face(\z) \not > 0$ for any face $\face$ adjacent to $\vertex \in \vertices_0$.

We can obtain an even more natural approach when we work with the nonlinear membrane energy $\W_\mem$.
Recall that in \autoref{sec:energy} we introduced the characterization \autoref{eq:manifoldViaEnergy} of $\manifold$ without explicit dependence on the triangle inequalities by exploiting the growth of $W_\mem$ for triangles with vanishing area.
This now readily fits into our modified line search approach.
In fact, if we compare our nonlinear energy to interior point methods \cite[Chapter 19]{NoWr06} we see that the logarithmic penalty in the energy takes the role of a barrier term which ensures that we stay in the admissible set $\trineqSet$.

Overall, we see that in both cases we can treat the inequality constraints in the line search and hence apply algorithms for equality-constrained optimization with a considerably lower number of constraints.
Note, that this approach can be adapted for trust region methods by limiting the size of the trust region appropriately.

\paragraph{Augmented Lagrange}
Next, we describe our approach to solving these equality-constrained problems based on the augmented Lagrange method.
First, let us briefly recall the Lagrangian formulation of our problem.
In fact, this means we seek for a saddle point of the Lagrangian 
	\begin{equation} \tag{Lag} \label{eq:genericLagrangian}
	\lagrange(\z, \mult) = \energy(\z) - \qint(\z) \cdot \mult
	\end{equation}
with $\z\in \R^{2\numE}$ and Lagrange multiplier $\mult \in \R^{3|\vertices_0|}$. 
The necessary condition for a saddle point $(\z,\mult) \in \R^{2\numE} \times \R^{3|\vertices_0|}$ is 
	\begin{align} \label{eq:lagrangeNecessary}
	D \lagrange(\z,\mult) = \left(D_\z \lagrange(\z,\mult),\, D_\mult \lagrange(\z,\mult)\right)^T = \left(D_\z \energy(\z) - D_\z \qint(\z) \cdot \mult,\, -\qint(\z)\right)^T = 0\,,
	\end{align}
where $D_\z$ and $D_\mult$ denote the Jacobian with respect to $\z$ and $\mult$, respectively.

Instead of directly applying Newton's method to this equation we consider the augmented Lagrange method \cite{He75, NoWr06}. 
It is a combination of the Lagrangian approach with the quadratic penalty method where we construct a series of unconstrained optimization problems in $\z$ to approximate the solution of \autoref{eq:genericOpt}.
For the sake of completeness, we briefly recall it here.
The augmented Lagrangian is defined by 
	\begin{equation} \label{eq:augmentedLagrangian}
	\lagrange(\z, \mult, \pen) = \energy(\z) - \qint(\z) \cdot \mult + \frac{\pen}{2}\, \lVert \qint(z) \rVert_2^2,
	\end{equation}
and a sequence $(\z^k, \mult^k, \pen^k)$ of approximate solutions, approximate Lagrangian multipliers, and penalty parameters is generated by alternating between minimizing $\lagrange(\,\cdot\, , \mult^k, \pen^k)$ to obtain $\z^{k+1}$ and computing updates to  $\mult^k$ and $\pen^k$.
Hereby, the penalty parameter $\pen$ is increased until we reach sufficient attainment of the equality constraints.
On the other hand, $\mult$ is updated by an increasingly accurate estimation of the correct multipliers $\mult^\ast$ solving \autoref{eq:lagrangeNecessary}.
This can be accomplished in various ways, one popular way which we choose to follow is to set $\mult^{k+1} = \mult^{k} - \pen^k\, \qint(z^{k+1})$.
Though we cannot expect the augmented Lagrange method to converge for arbitrary initial data, 
under reasonable assumptions, one can prove that the sequence $\mult^k$ obtained this way converges to $\mult^\ast$, which significantly improves convergences compared to the quadratic penalty method, see for example \cite{NoWr06}.
We want to remark that though our problem \autoref{eq:genericOpt} involves strict inequality constraints, the local convergence theory for the augmented Lagrange method given in \cite[Chapter 17]{NoWr06} applies to our problem. 
The triangle inequality constraints define an open set and thus \autoref{eq:genericOpt} can be seen as an equality-constrained problem over an open set.
As \cite[Theorem 17.5 \& 17.6]{NoWr06} are only concerned with local minimizers and provide local results, they still hold if the problem is only defined on an open set after possibly modifying constants describing local neighborhoods.

An explicit algorithmic description of the method with all involved parameters and derivatives will be provided in the appendix.

\paragraph{Unconstrained Optimization}
Using the augmented Lagrange method leads to a series of unconstrained optimization problems.
They are typically non-convex, \ie we encounter indefinite Hessians $D^2_\z \lagrange$ of the Lagrangian.
This means that a simple Newton's method with line search might not be an efficient and robust approach as we are not guaranteed to obtain a descent direction.
To rectify this, we choose a simple adaption suggested in \cite[Section 3.4]{NoWr06}.
First, we determine a shift $\shift^k$ such that the matrix $ D^2_\z \lagrange(\z^k, \mult^k, \pen^k) + \shift^k \Id$ is positive definite.
This achieved by starting with an initial estimate and then increasing $\shift^k$ until a Cholesky decomposition succeeds.
Then, a descent direction is obtained by solving the linear system
\begin{equation}
	\label{eq:descentDir}
	\left(D^2_\z \lagrange(\z^k, \mult^k, \pen^k) + \shift^k \Id\right)d^k = -D_\z \lagrange(\z^k, \mult^k, \pen^k).
\end{equation}
Along this direction we perform an Armijo-type backtracking line search.
Note again, that the local convergence theory for Newton-type methods is still valid even though we minimize over an open set defined by the strict triangle inequalities, \cf \cite[Chapter 1]{Be99}.
In some instances, we could speed-up the minimization by first performing a small number of iterations with a BFGS approximation of the Hessian. 

To compute the descent direction as above, we need the gradient and the Hessian of our constraint functionals $\qint$. 
We already evaluated $D_\z \qint \in \R^{3|\vertices_0|, 2\numE}$ in \autoref{eq:intFirstDeriv} and compute for the Hessian of $\qint$
\[
D_\z^2 \qint \cdot \mult = \left(\sum_{\vertex \in \vertices_0}  \partial_{\z_l} \partial_{\z_k} \qint_\vertex \cdot \mult_\vertex \right)_{l,k = 1,\ldots, 2\numE}
\]
the components as
\begin{align*}
\partial_{z_l} \partial_{z_k} \qint_\vertex(z) &= \mathrm{vec} \left(\sum_{j=0}^{n-1} q_{01}(\z) \ldots  \partial_{z_l} \partial_{z_k} q_{j,j+1}(\z)  \ldots  q_{n-1,0}(\z)\right)\\ 
&\phantom{=} + \mathrm{vec} \left(\sum_{\substack{i,j=0\\i\neq j}}^{n-1} q_{01}(\z) \ldots \partial_{z_l} q_{i,i+1}(\z)  \ldots  \partial_{z_k} q_{j,j+1}(\z)   \ldots  q_{n-1,0}(\z)\right),
\end{align*}
which can also be evaluated with $O(\nv)$ cost. 
We provide further details on the Hessian computation in the supplementary material.

\section{Reconstruction of an immersion}
\label{sec:recon}
In the preceding sections, we discussed the geometry as well as constrained optimization problems on the NRIC manifold $\manifold$, \ie in terms of edge lengths and dihedral angles.
The remaining task is to reconstruct for given $\z \in \manifold$ an immersion $\pos \in \vertexSpace$ of the simplicial surface in $\R^3$ with $\z=\projZ(\pos)$. 
Beyond the computation of vertex coordinates for  $\z\in\manifold$, one frequently asks for an approximate immersion $\pos\in\vertexSpace$ for $\z \not\in \manifold$ such that $\projZ(\pos) \approx \z$.
Indeed, the computation of just approximate immersions is required in case of  
\begin{itemize}
	\item modeling of deformations energies in terms of dihedral angles and edge lengths, \ie using the linear embedding space $\R^{2\numE}$ instead of $\manifold$,
	\item using a high tolerance for the fulfillment of the constraints in the augmented Lagrange or a penalty method,
	\item coordinates $\z$ which are only approximately computed numerically. 
\end{itemize}
Thus,  we ask for a reconstruction map \(\rec \colon \R^{2\numE} \to \vertexSpace\), such that \(\rec\vert_\manifold\) is the right inverse of \(\projZ\) with
\(\projZ \circ \rec = \Id_ \manifold\), where \(\Id_ \manifold\) is the identity on the NRIC manifold. 
Let us emphasize that by the rigid body motion invariance of our NRIC approach,  we obtain \(\rec \circ \projZ (\pos) =  Q \pos\), where \(\pos \in \vertexSpace\) and \(Q \in \mathit{SE}(3)\) is some rigid body motion acting on the immersion.

\paragraph{Variational approach}
For some given $\z \in \R^{2\numE}$, where not necessarily $\z \in \manifold$, we are looking for the nodal positions $\pos \in \vertexSpace$, such that the resulting $\projZ(\pos) \in \manifold$ is as close as possible to $\z$. 
Fr\"ohlich and Botsch \cite{FrBo11} have used a least squares functional to build a \emph{variational reconstruction}, \ie they compute
\begin{equation}\label{eq:genericVarRecon}
\argmin_{\pos \in \R^{3\numV}} \, \W\left[\z,\,\projZ(\pos)\right] 
\end{equation}
with $\W\left[\z,\tilde \z\right]$ describing the proximity of  $\z$ and $\tilde \z$.
Note that the solution is only unique up to a rigid deformation. 
A simple example of a quadratic functional $\W$ is given by \eqref{eq:quadEnergy} as it was used by Fr\"ohlich and Botsch. 
They proposed a Gau\ss-Newton method \cite[Section 10.3]{NoWr06} to solve \eqref{eq:genericVarRecon}, however, for general  $\z\in \R^{2\numE}$, one still has to solve a high-dimensional and nonlinear optimization problem in $\R^{3\numV}$.
If \(\z \in \manifold\) and the initialization of the Gau\ss-Newton method is close to the solution, it usually converges in only a few iterations.
However, if $\z$ is far away from \(\manifold\) and the initialization is poor, artifacts may occur. 

\paragraph{Constructive approach}
For $\z \in \manifold$ a \emph{constructive reconstruction} of the immersion $\pos\in\vertexSpace$ can be derived by means of frames and transition rotations,  as they were used to define the integrability conditions (\cf Section \ref{sec:background}). 
This method was introduced by \citet{LiSoLe05} and further elaborated in \cite{WaLiTo12}.
Before we investigate a combination of the constructive and the variational approach for $\z\notin\manifold$, let us briefly review the constructive reconstruction. 
Assume we are given an admissible target \(\z = (\len,\dih) \in\manifold\).
Since the reconstruction from lengths and angles is only defined up to rigid body motions, we further assume that we are given the position of one vertex and the orientation of an adjacent triangle \(\face_0\) in the form of a frame \(\frame_0\). 
If $\face\in\faces$ is a neighboring triangle of \(\face_0\), one can infer the induced transition rotations $R_{0\face}$ from $\z$ and thus determine  $\frame_\face = \frame_0 R_{0\face}$. Repeating this iteratively, one can construct frames for all faces.
This algorithm is indeed well-defined on simply connected triangulations due to the integrability constraints, \ie  if there are two paths connecting a triangle \(\face\) to \(\face_0\), then the frames constructed along the two paths coincide.
Given frames for all faces and hence the orientation of all triangles, one can finally reconstruct the nodal positions.

\paragraph{Adaptive spanning trees}
Next, we take into account a violation of the integrability condition \eqref{eq:discreteIntegrMap} for $\z \not\in \manifold$ and ask for a reconstruction of an approximate immersion. 
Let us remark that \citet{WaLiTo12} handle non-admissible targets $\z \notin \manifold$ when modeling surfaces via a modification of  $\len$ and $\theta$.
They study a least-square type functional and relax in a least square sense the identity $\frame_\face = \frame_0 R_{0\face}$ as well as \eqref{eq:edgeVectors}.

\begin{figure}[ht]
	\centering
	\includegraphics[width=\columnwidth]{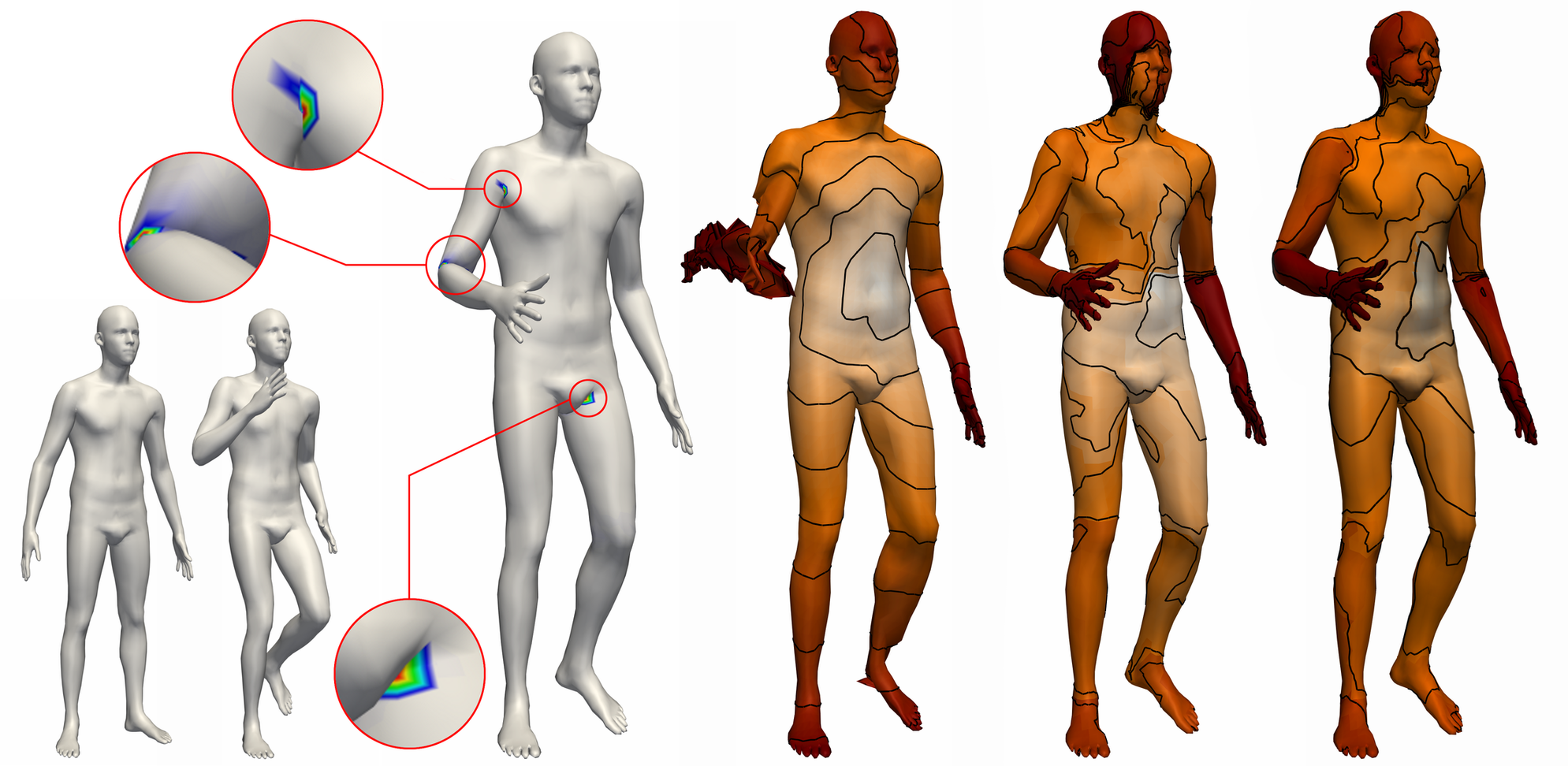}
	$\pos_1$ \hspace*{12mm} $\pos_2$ \hspace*{5.5cm} (BFS) \hspace*{17mm} (MST) \hspace*{17mm} (SPT)
	\caption{Left: Input shapes $\pos_1$ and $\pos_2$ (taken from \cite{PoRoMaBl15}) and reconstruction from linear average $(\projZ(\pos_1) + \projZ(\pos_2))/2 \notin \manifold$ with the local violations of the integrability condition as color map. Note that violations are highly concentrated, \eg in the armpits.
		Rightmost shapes: reconstruction using breadth-first search (BFS), minimal spanning tree (MST) and shortest path tree (SPT). The triangulation is color-coded with respect to the order of traversal. 
		See video for an animation of the reconstruction order.}
	\label{fig:dynaReconstruction}
\end{figure}

The direct frame-based reconstruction with a spanning tree of the dual graph built by breadth-first search is very sensitive to violations of the integrability.
In fact, the errors occurring when walking over such a violation propagate to all following frames and are even amplified, \cf \autoref{fig:dynaReconstruction}.
In addition, reconstructing nodal positions of a face along two different paths connecting it to the initial face \(\face_0\), where at least one is passing a zone of violated integrability conditions, leads to substantially different results and thus visual artifacts.
However, the regions of violation appear frequently to be highly localized in practice,  \cf \autoref{fig:dynaReconstruction}.
Thus, we build a spanning tree which traverses faces with violation of the integrability condition as late as possible in the mesh traversal for the reconstruction.
To this end, we consider the dual graph of $\simplicial$ with weights based on the integrability condition.
Each dual edge corresponds to a primal edge \(\edge = (\vertex,\vertex ') \in \edges\) and we can assign to this dual edge a scalar weight reflecting the lack of integrability \(w_\edge\) by averaging the violation of integrability at the two adjacent vertices \(\vertex\) and \(\vertex'\):
\begin{equation}\label{eq:edgeWeights}
w_\edge := \frac{\rvert\tr\,  \dint_{\vertex}(\z) - 3\lvert + \rvert\tr\,  \dint_{\vertex'}(\z) - 3\lvert}{2}\,,
\end{equation}
where $\dint_{\vertex}(\z)$ is the matrix-valued map defined in \eqref{eq:discreteIntegrMap}.
Note that $\dint_{\vertex}(\z) \in \SO(3)$ and that $\tr Q = 3 \Leftrightarrow Q= \Id$ for $Q \in \SO(3)$. 
Now, the weights \eqref{eq:edgeWeights} are used to build a spanning tree adapted to the problem. 
The first variant is to construct a minimal spanning tree (MST) of the dual graph, which is built such that the sum of all edge weights in the tree is minimal.
Such a minimal spanning tree can be computed via Prim's algorithm and provides a way to traverse the dual graph while avoiding unnecessarily large violations of the integrability.
Another variant is to construct a shortest path tree (SPT), which is built such that the path distance from the root to any other vertex in the tree is the shortest in the whole dual graph.
This can be achieved by Dijkstra's algorithm and provides a way to traverse the dual graph such that for each face the sum of integrability violation along the dual path used for its reconstruction is minimal. 
We compare both novel variants against the original breadth-first search (BFS) in \autoref{fig:dynaReconstruction}. 
A pseudo code of the entire algorithm is given in the appendix.
Formally, the algorithm --- using either (MST) or (SPT) --- is only defined for $\z\in\manifold$. 
In particular, the triangle inequality is assumed to be defined.
However, the algorithm can easily be generalized for $\z \in \R^{2\numE}$ with $\trineq_\face(\z) \leq 0$ for some face $\face$ by setting the interior angles of $\face$ to zero. 
By our definition of \eqref{eq:discreteIntegrMap} and edge weights \eqref{eq:edgeWeights} those triangles will be automatically considered as late as possible in the adaptive algorithm.

\paragraph{Preassembled tree}
The runtime of the tree-based reconstruction algorithm is dominated by the cost for the construction of the spanning tree.
Thus, if one aims at reconstructing numerous immersions of a discrete surface with the same connectivity and a very high resolution (\ie many vertices) it would be desirable to use a preassembled spanning tree.
Of course, this preassembled tree has to be reasonable for a large set of lengths and angles. 
If we are given samples $z_1, \ldots, \z_n \in \manifold$ and corresponding edge weights $w^1, \ldots, w^n \in \R^{\numE}$, we simply set $ w_\edge = \max_i w^i_\edge$ for all $\edge \in \edges$ and construct a spanning tree based on these weights.

\begin{figure}[t]
	\includegraphics[width=\columnwidth]{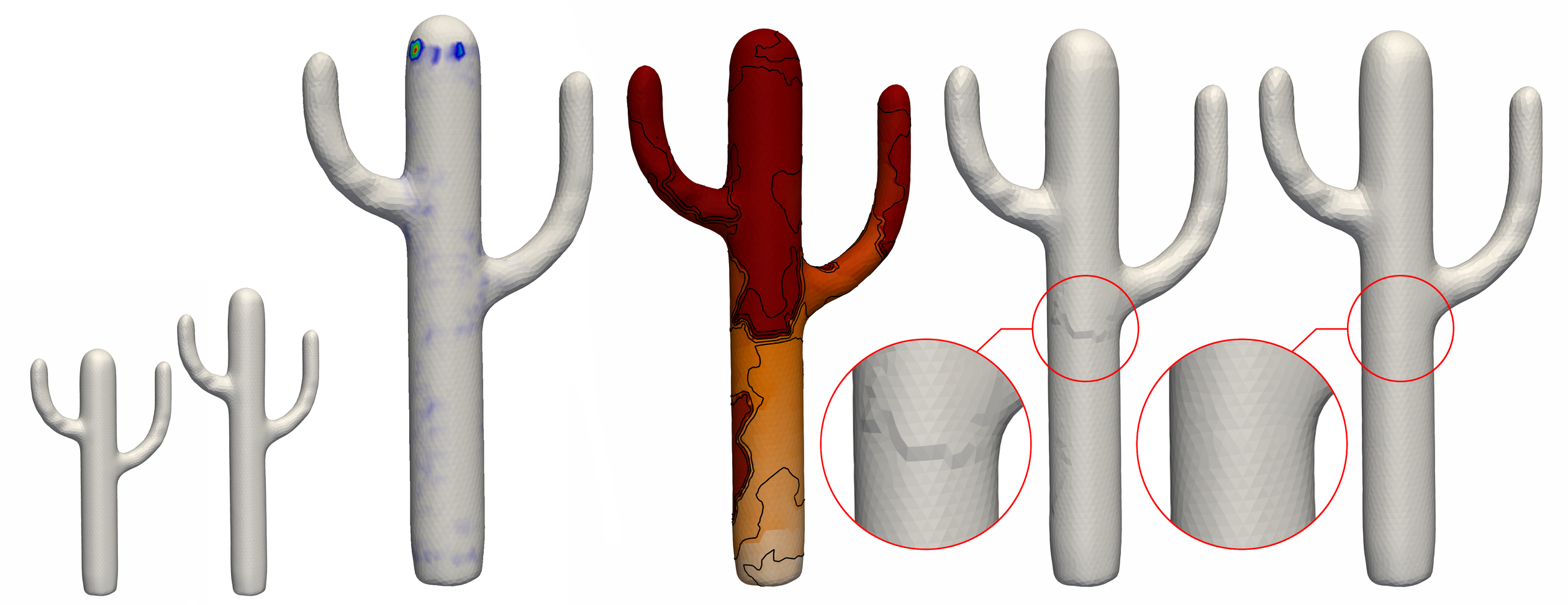}
	\caption{Left: Input shapes $\pos_1$ and $\pos_2$ (small) and reconstruction from linear average $(\projZ(\pos_1) + \projZ(\pos_2))/2 \notin \manifold$ with local violations of the integrability condition as color map. 
		Rightmost shapes: reconstruction order using a minimal spanning tree (order as colormap), visual artefacts along the body and final solution after one step of a Gau\ss-Newton smoothing.}
	\label{fig:cactusReconstruction}
\end{figure}

\paragraph{Hybrid approach}
Just applying our constructive reconstruction algorithm works very well for $\z\notin\manifold$ as long as the violations are localized as in \autoref{fig:dynaReconstruction}. 
However, we observe imperfect results when the violations are distributed over larger areas, \cf \autoref{fig:cactusReconstruction}. 
In this case, we suggest a hybrid method combining our robust constructive reconstruction and as a post processing the variational reconstruction. 
In detail, we make use of the (still imperfect) output of our constructive reconstruction to initialize the variational reconstruction as in \eqref{eq:genericVarRecon}. 
Typically, a single Gau\ss-Newton step is sufficient to smooth the result adequately (\cf \autoref{fig:cactusReconstruction}).  

\section{Numerical experiments and comparisons}\label{sec:results}
In this section, we study qualitative and quantitative properties of the NRIC tools and demonstrate that in particular for modeling with near isometric deformations the NRIC manifold outperforms established methods that consider nodal positions as primal degrees of freedom. 
To this end, we pick up the generic variational problem \autoref{eq:genericOpt} introduced in \autoref{sec:varProblems} together with the proposed augmented Lagrange method. 
In the following, we discuss different objective functionals $\energy$ in \autoref{eq:augmentedLagrangian} and depending on the application additional constraints. 
Note, however, that the \emph{constraint functional} $\qint$ in \autoref{eq:augmentedLagrangian}, which describes the NRIC manifold implicitly via \autoref{eq:nricDefinitionImplicit}, remains unchanged. 

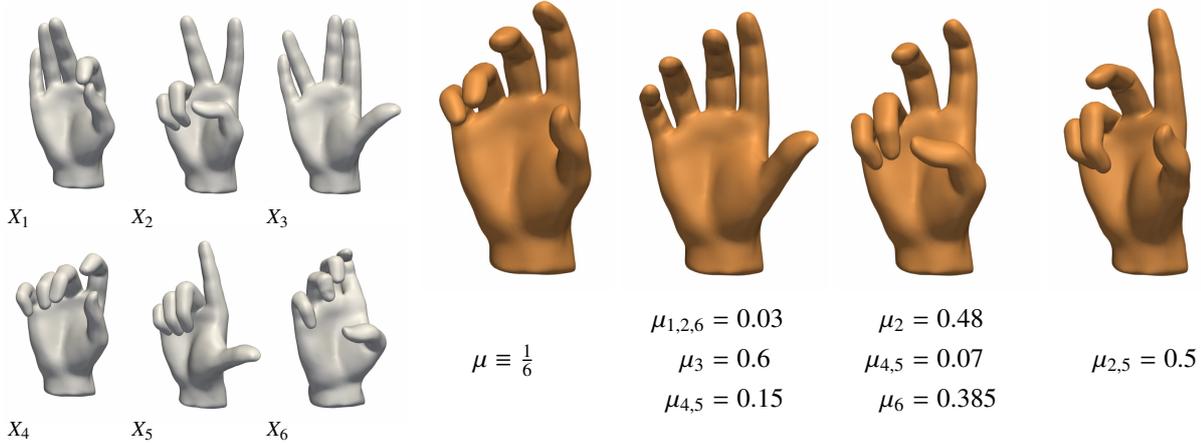
\begin{figure}[ht]
	\begin{minipage}[c]{0.3\columnwidth}
		\setlength{\unitlength}{.3\linewidth}%
		\begin{picture}(3,4)
		\put(-0.2,2){ \includegraphics[width=1.2\unitlength]{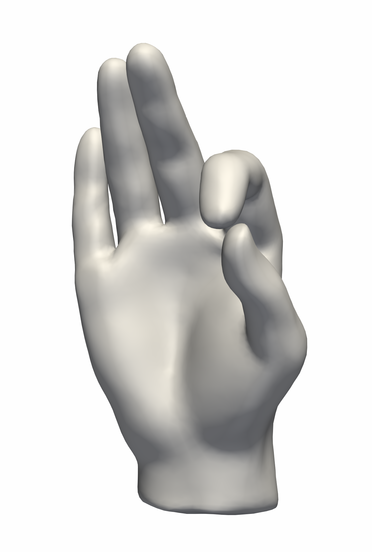} }
		\put( 1.0,2){ \includegraphics[width=1.2\unitlength]{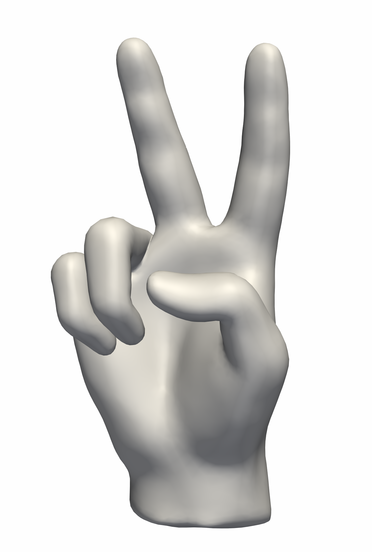} }
		\put( 2.2,2){ \includegraphics[width=1.2\unitlength]{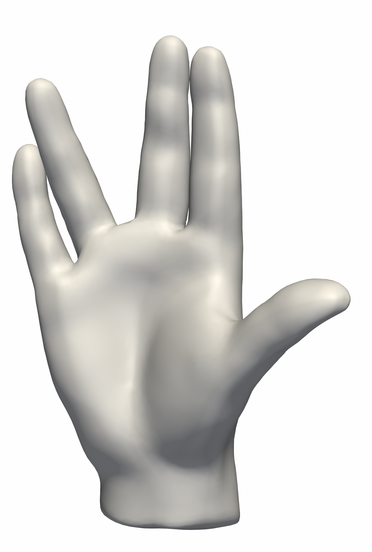} }
		\put(-0.2,0){ \includegraphics[width=1.2\unitlength]{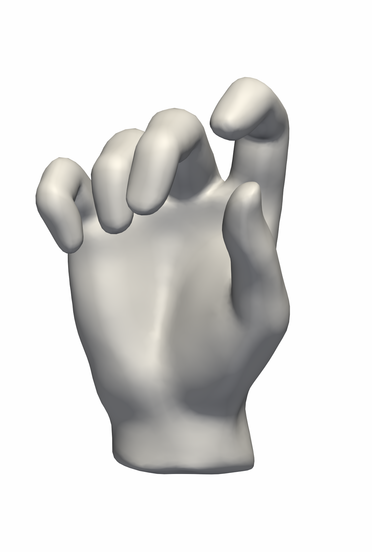} }
		\put( 1.0,0){ \includegraphics[width=1.2\unitlength]{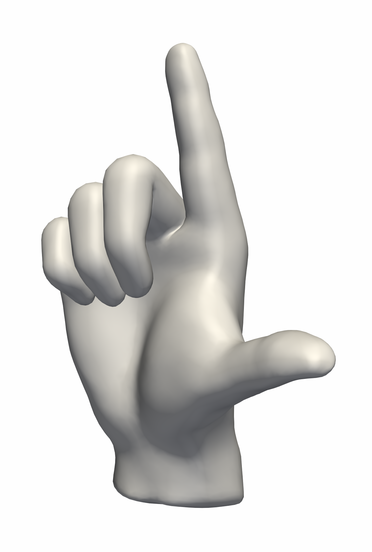} }
		\put( 2.2,0){ \includegraphics[width=1.2\unitlength]{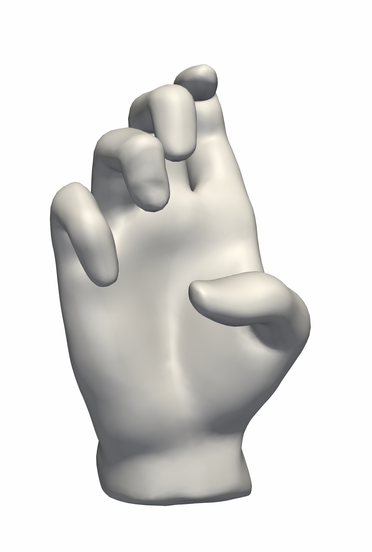} }
		\put(-0.1,1.8){\footnotesize $\pos_1$ }
		\put( 1.0,1.8){\footnotesize $\pos_2$ }
		\put( 2.2,1.8){\footnotesize $\pos_3$ }
		\put(-0.1,-0.1){\footnotesize $\pos_4$ }
		\put( 1.0,-0.1){\footnotesize $\pos_5$ }
		\put( 2.2,-0.1){\footnotesize $\pos_6$ }
		\end{picture}
	\end{minipage}
	\begin{minipage}[c]{0.69\columnwidth}
		\setlength{\fboxsep}{0pt}%
		\setlength{\fboxrule}{1pt}%
		\hspace{0.5em}
		\includegraphics[width=0.24\linewidth-0.5em,trim={0 1em 2cm 0},clip]{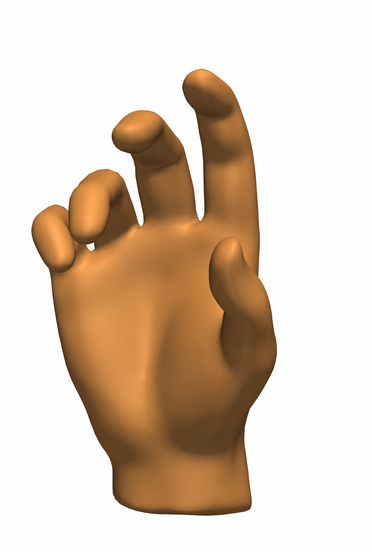}
		\includegraphics[width=0.24\linewidth,trim={0 1em 0 0},clip]{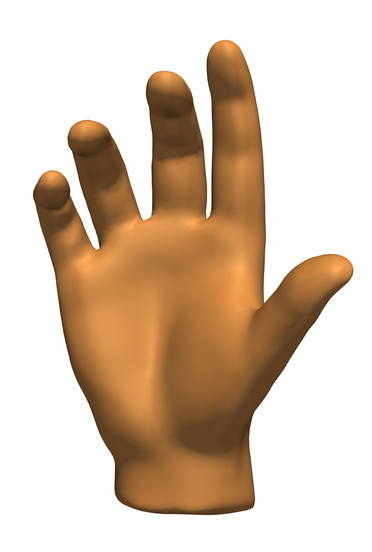}
		\includegraphics[width=0.24\linewidth,trim={0 1em 0 0},clip]{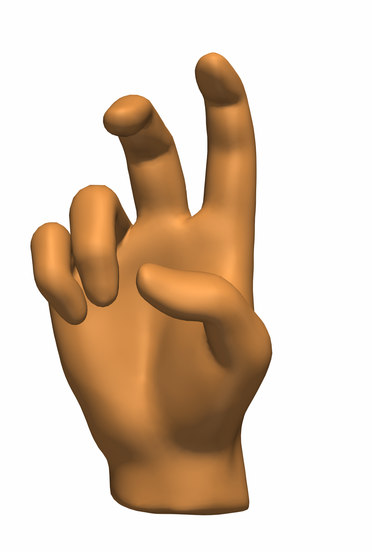}
		\includegraphics[width=0.24\linewidth,trim={0 1em 0 0},clip]{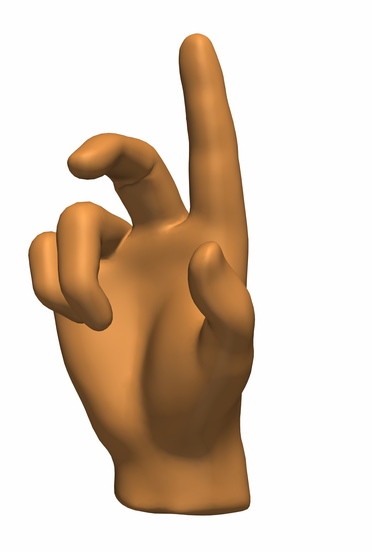}\\[0.1em]
		\parbox{0.24\linewidth}{\vspace{-1em}\begin{align*} \mu \equiv \tfrac16 \end{align*}}
		\parbox{0.24\linewidth}{\vspace{-1em}\begin{align*} \mu_{1,2,6} &= 0.03 \\ \mu_3 &= 0.6 \\ \mu_{4,5} &= 0.15 \end{align*}}
		\parbox{0.24\linewidth}{\vspace{-1em}\begin{align*} \mu_{2} &= 0.48 \\ \mu_{4,5} &= 0.07 \\ \mu_{6} &= 0.385 \end{align*}}
		\parbox{0.24\linewidth}{\vspace{-1em}\begin{align*} \mu_{2,5} = 0.5 \end{align*}}
		
	\end{minipage}
	\caption{Reconstruction of nodal positions from elastic averages of six hand poses (grey) with different (convex) weights $\mu\in\R^6$ computed as minimizer of \eqref{eq:elasticAverageFunctional} on the NRIC manifold. }
	\label{fig:handsAverage}
\end{figure}

\paragraph*{Elastic averages}
Let $\pos_1, \ldots, \pos_n \in \vertexSpace$ be a set of example shapes (sharing the same connectivity). 
Frequently, one is interested in a mean or average shape, \cf \cite{TyScSe15}. 
Given an elastic deformation energy, a so-called weighted elastic average is defined to be the minimizer of a weighted sum of elastic energies for deformations from the input shapes to the free shape.
This can be translated directly to our NRIC manifold, \ie for a given elastic deformation energy $\W$ on $\manifold$ and convex weights $\mu \in \R^n$ we define the weighted elastic NRIC average as a solution of \autoref{eq:genericOpt} with 
\begin{equation} \label{eq:elasticAverageFunctional}
\energy(\z) = \sum_{i = 1}^n \mu_i \, \W[ \projZ(\pos_i), \z ]\, .
\end{equation}
In \autoref{fig:handsAverage}, we show (the reconstructions of) weighted elastic NRIC averages for a set of six hand shapes and different weights $\mu_1, \ldots, \mu_6$. 
Here, we have used the nonlinear deformation energy \eqref{eq:nonlinearEnergy} in \eqref{eq:elasticAverageFunctional}.

\begin{figure}[ht]
	\centering
	\includegraphics[width=0.48\columnwidth]{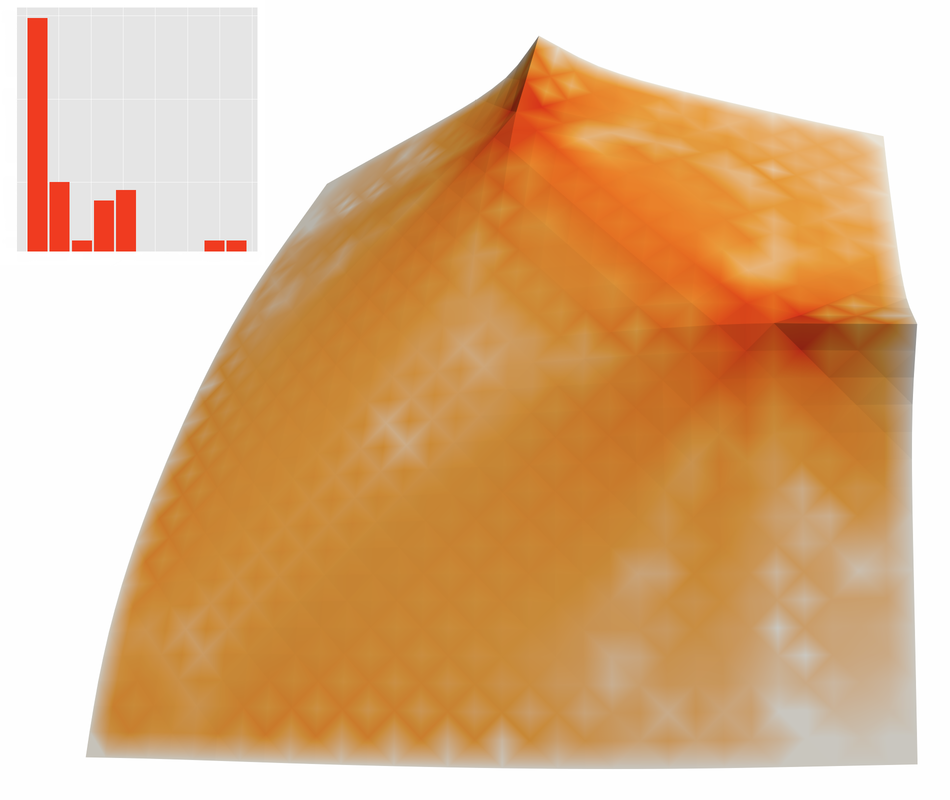}
	\includegraphics[width=0.48\columnwidth]{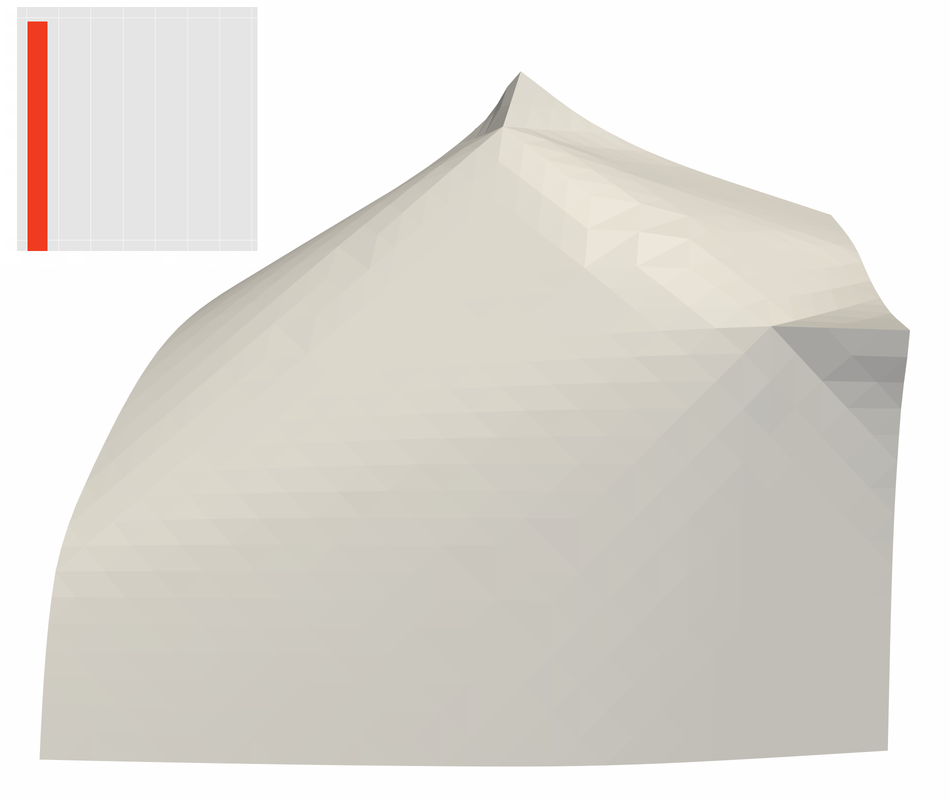}
	\caption{Paper folding with local constraints for dihedral angle: simulation in vertex space (left) leads to infinitesimal isometry violations whereas the result in NRIC is completely isometric (right). 
		The absolute value of discrete Gau\ss \, curvature (as angle defect) is shown using the color map $0\hspace{1mm}$\protect\resizebox{.1\linewidth}{!}{\protect\includegraphics{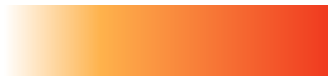}}$\hspace{1mm}0.03$, which is zero everywhere on the right. 
		Furthermore, the corresponding histograms are plotted aside the surfaces.}
	
	\label{fig:paperFoldingCompareGauss}
\end{figure}

\paragraph*{Isometric deformations via additional constraints} 
Interesting applications can be described by considering \autoref{eq:genericOpt} along with the simple objective $\energy(\z) = \W[\z^\ast,z]$ but with additional, simple coordinate constraints. 
For example, in \autoref{fig:strangulation} we have seen experiments where we posed lengths constraints $\len^{\vphantom{ast}}_{i} = \len^\ast_i$ for $i\in I$ on the coordinates $\z = (\len,\dih)$ for some index set $I \subset \edges$ and prescribed target lengths $\len^\ast$. 
Similarly, we obtain an elegant way to simulate the \emph{isometric} folding of a (flat) sheet of paper given in NRIC as $\z^\ast = (\len^\ast, \dih^\ast)$ where $\dih^\ast = 0$.  
To this end, we pose the length constraints $\len^{\vphantom{ast}}_\edge = \len^\ast_\edge$ for all $\edge \in \edges$ along with $\dih_{i} = const \neq 0$ if $i\in I$ for some index set $I \subset \edges$. 
Note, that under these length constraints the nonlinear and quadratic energy approach agree if we compute the weights in \autoref{eq:quadEnergy} from the reference $\z^\ast$.
For example, in \autoref{fig:paperFoldingCompareGauss} we impose the constraint $\dih_i = \pi/2$ for the edges on two short line segments on two neighboring sides of the sheet. 
Since all edge lengths are fixed and all other dihedral angles are degrees of freedom for the minimization of \autoref{eq:nonlinearEnergy} on $\manifold$, we obtain a perfect isometric deformation as indicated by the vanishing discrete Gau\ss \, curvature (\autoref{fig:paperFoldingCompareGauss}, right). 
In comparison, vertex-based methods as \cite{GrHiDeSc03} or \cite{HeRuWa12} do not achieve a perfect isometry---even when computed with a very high membrane stiffness (\autoref{fig:paperFoldingCompareGauss}, left). 
For the optimization in nodal positions, we used the energy $\pos \mapsto  \W_{nl}[\z^\ast, \projZ(\pos) ]$ with a shell thickness parameter $\delta=10^{-3}$. 
In fact, further reducing $\delta$ one observes numerical instabilities.
This is due to the fact that isometric deformations induce bending distortions only but optimizing bending energies in terms of nodal positions is a highly nonlinear singular perturbation problem that quickly triggers numerical issues. 
Conversely, the corresponding bending energy in NRIC is quadratic.

\begin{figure}[ht]
	\centering
	\includegraphics[width=0.32\columnwidth]{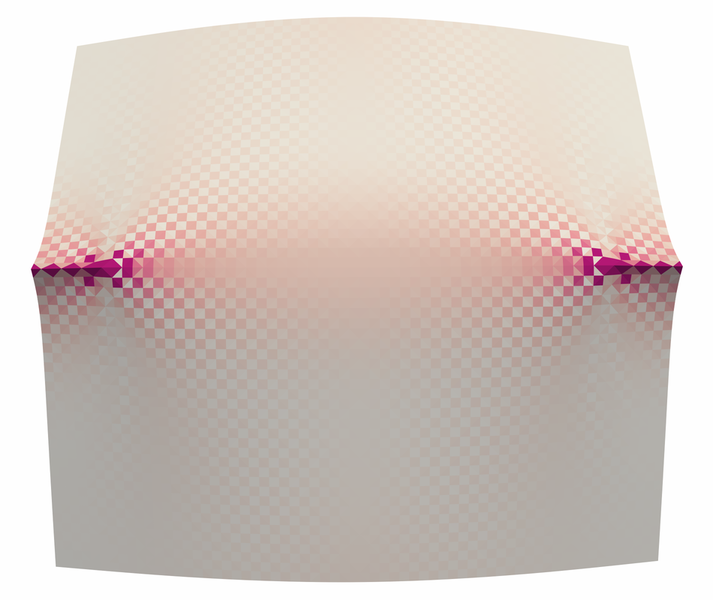}
	\includegraphics[width=0.32\columnwidth]{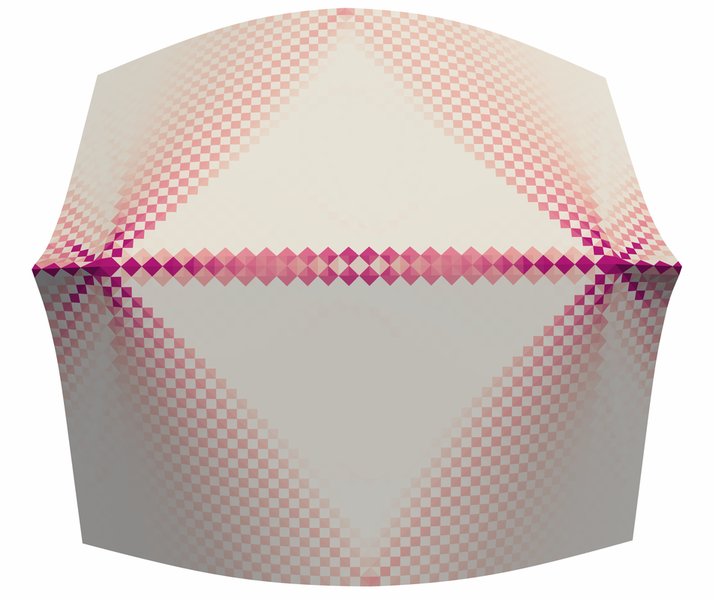}
	\includegraphics[width=0.32\columnwidth]{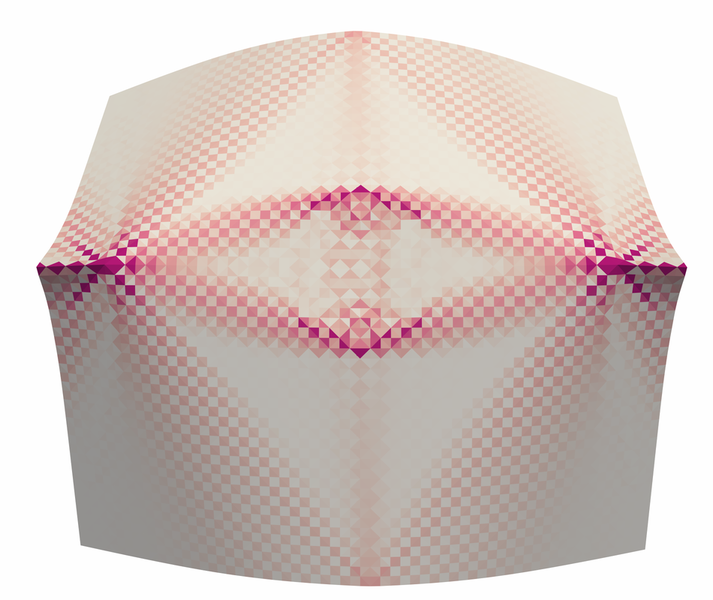}\\
	\includegraphics[width=0.32\columnwidth]{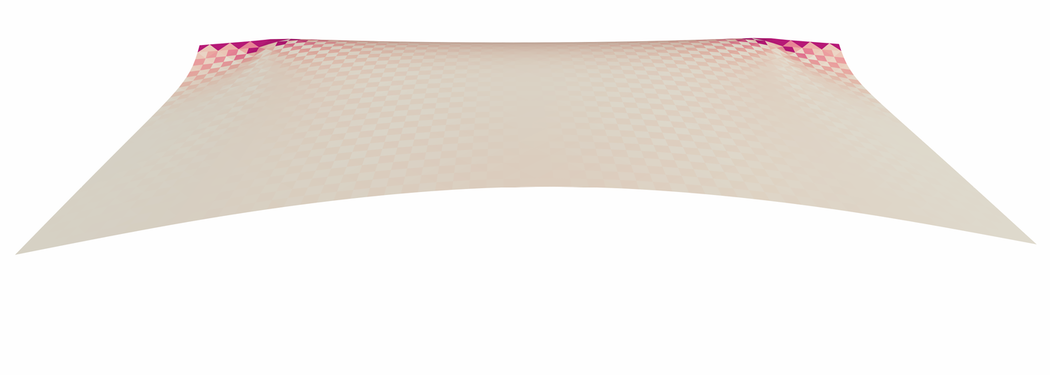}
	\includegraphics[width=0.32\columnwidth]{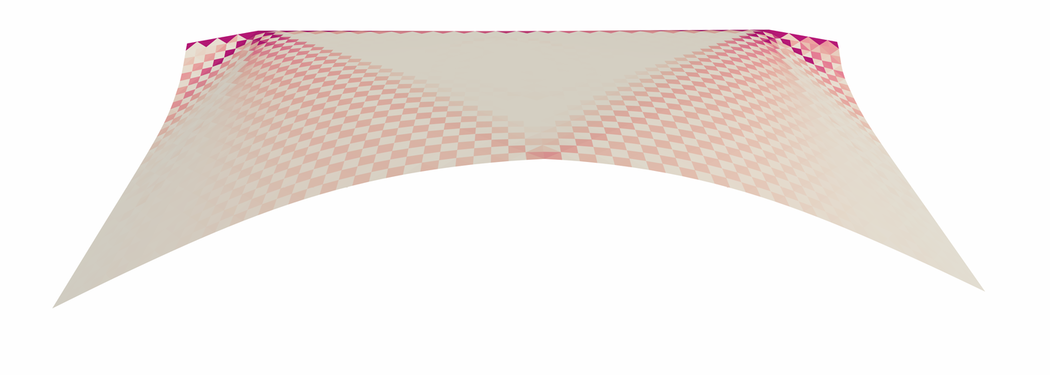}
	\includegraphics[width=0.32\columnwidth]{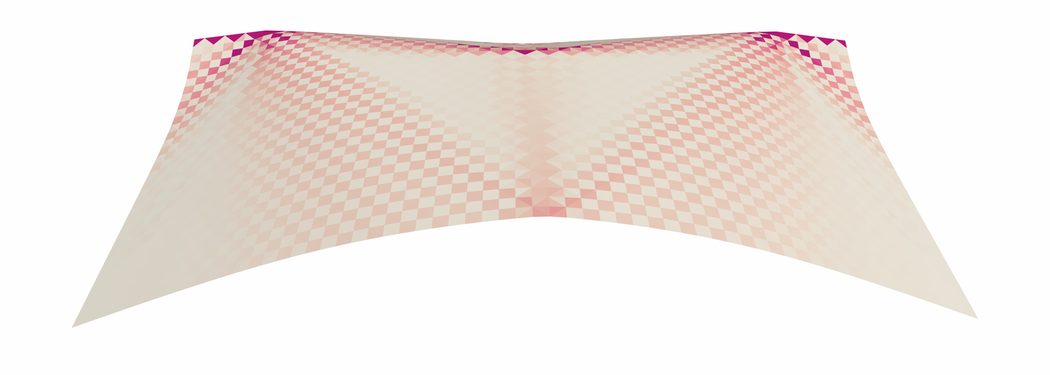}\\
	\parbox{0.32\columnwidth}{\hfill}
	\parbox{0.32\columnwidth}{\centering $\energy(\z) \approx 159.6$}
	\parbox{0.32\columnwidth}{\centering $\energy(\z) \approx 150.3$}
	\caption{Once more paper folding with local constraints for dihedral angle: simulation in vertex space (left) leads to infinitesimal isometry violations whereas the results in NRIC are completely isometric (middle, right). 
		The result in the middle shows a local minimum obtained by the augmented Lagrange method when increasing the penalty parameter $\pen$ (too) aggressively exhibiting higher deformation energy (shown below) than the right result where the penalty was increased more conservatively.
 		Triangle-averaged mean curvature is shown as color map  $0\hspace{1mm}$\protect\resizebox{.1\linewidth}{!}{\protect\includegraphics{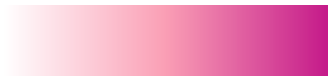}}$\hspace{1mm}\geq0.005$, flat triangular regions can only be seen in the NRIC simulations. }
	\label{fig:paperFoldingCompare}
\end{figure}

Besides vanishing Gau\ss \, curvature, pure isometric deformations exhibit further characteristics, as illustrated in \autoref{fig:paperFoldingCompare}. 
In this example, we have a very similar setup as in \autoref{fig:paperFoldingCompareGauss} but we pose the angle constraints on two opposite sides. 
First, let us point out that we observed convergence of the augmented Lagrange method described above to different local minima when using different parameters for the increase of the penalty parameter $\pen$.
We show two different local minima in \autoref{fig:paperFoldingCompare} where we obtained the lowest energy value when increasing $\pen$ conservatively (shown on the right).
Now, since the NRIC results are perfectly isometric and rather smooth deformations of the flat sheet one can indeed observe effects predicted analytically by the Hartman-Nirenberg theorem \cite{HaNi59, Ho11}. 
Loosely speaking, isometric deformations of a flat sheet can locally be described either as flat patches or segments of straight lines (rulings) going to the boundary. 
In the middle and right columns of \autoref{fig:paperFoldingCompare} one can easily identify flat triangular regions as well as a cone-like bundle of straight lines propagating towards the boundary.
These structures are not reflected by the vertex-based numerical minimizer already discussed above (\autoref{fig:paperFoldingCompare}, far left).

\paragraph*{Time-discrete geodesics}
So far we have only considered static examples where a single shape was optimized subject to external forces or boundary conditions. 
However, one can easily generalize \autoref{eq:genericOpt} to optimize for multiple shapes simultaneously, for instance, to simulate a kinematic behavior. 
We focus on the computation of time discretized \emph{geodesics} in the NRIC manifold here. 
On the manifold  $\manifold$ with metric $\metric$ defined in \autoref{eq:metric} a geodesic connecting end points $\curve_A$ and $\curve_B$ in $\manifold$ is the curve $\curve \colon [0,1] \to \manifold$ minimizing the  
\emph{path energy} $\int_0^1 \metric_{\curve(t)}(\dot\curve(t),\dot\curve(t))\d t$ subject to  $\curve(0)=\curve_A$ and $\curve(1) = \curve_B$. 
In particular, the minimizer $(\curve(t))_{0\leq t\leq 1}$ obeys the constant speed property $\metric_\curve(\dot\curve(t), \dot\curve(t)) = const$. 
\citet{HeRuWa12} introduced the concept of \emph{time-discrete geodesics} (in a vertex-based approach) as a variational approximation of continuous geodesics. 
For $K\in\N$, they consider a finite sequence $\curve_0, \ldots, \curve_K$ in $\manifold$ with $\curve_0=\curve_A$ and $\curve_K=\curve_B$ and define the \emph{time-discrete path energy}
\begin{align}\label{eq:discretePathEnergy}
E[\curve_0, \ldots, \curve_{K}] = K \, \sum_{k=1}^K \W[\curve_{k-1}, \curve_k]\, ,
\end{align}
where $\W$ is assumed to be a local approximation of the squared Riemannian distance and $\curve_k \approx \curve(k/K)$.
Minimizers $(\curve_0, \ldots, \curve_K)$ of \eqref{eq:discretePathEnergy} for fixed end points $\curve_0$ and $\curve_K$ are said to be time-discrete geodesics. 
In particular, they obey a discrete constant speed property, \ie there is a uniform energy distribution $\W[\curve_{k-1}, \curve_k] \approx const$ along the curve. 

\begin{figure}[ht]
	\centering
	\scalebox{1}{
		\includegraphics[width=0.08\textwidth]{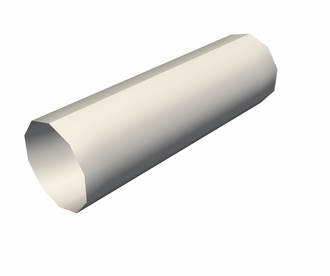}
		\includegraphics[width=0.08\textwidth]{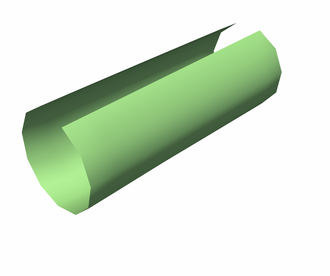}
		\includegraphics[width=0.08\textwidth]{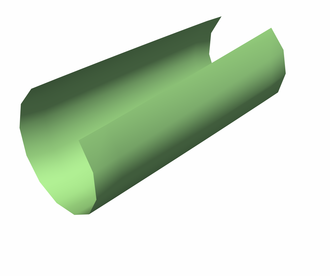}
		\includegraphics[width=0.08\textwidth]{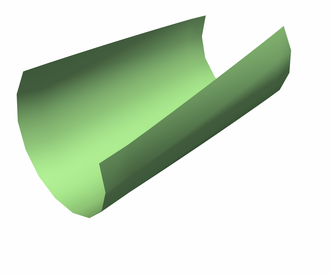}
		\includegraphics[width=0.08\textwidth]{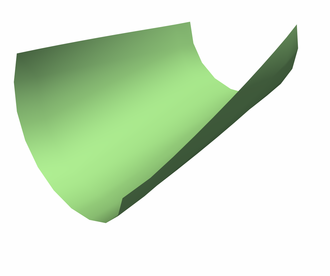}
		\includegraphics[width=0.08\textwidth]{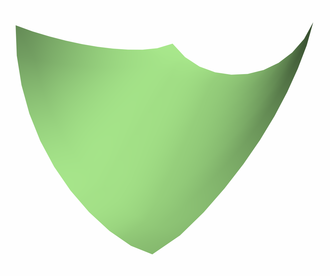}
		\includegraphics[width=0.08\textwidth]{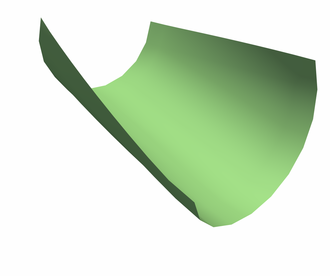}
		\includegraphics[width=0.08\textwidth]{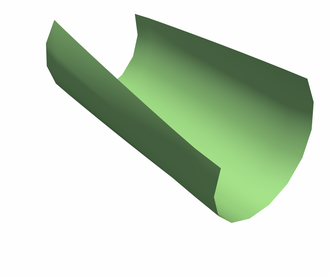}
		\includegraphics[width=0.08\textwidth]{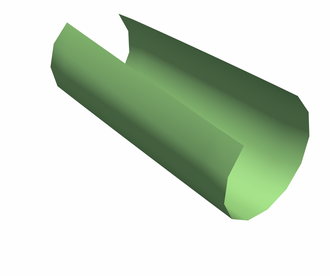}
		\includegraphics[width=0.08\textwidth]{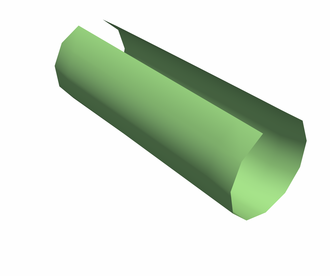}
		\includegraphics[width=0.08\textwidth]{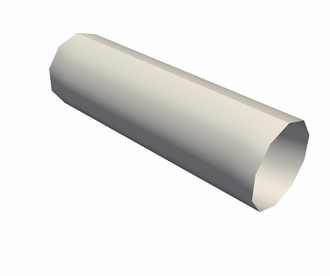}
	}\\[0.5em]
	\scalebox{1}{
		\begin{tikzpicture} 
		\begin{axis}[
		ymajorticks=false,
		axis y line=none,
		axis x line*=bottom,
		enlargelimits=0.05,
		enlarge x limits=0.001,
		ybar interval=0.8, 
		ymin=0.0015,
		ymax=0.023,
		width=0.95\textwidth,
		height=2cm 
		]
		\addplot[fill=myGreen!80!gray] coordinates {(0,0.001898) (1,0.001958) (2,0.002096) (3,0.002522) (4,0.007482) (5,0.022749) (6,0.002565) (7,0.002092) (8,0.001941) (9,0.001876) (10,0.001876)}; 
		\end{axis} 
		\end{tikzpicture}
	}\\[0.5em]
	\scalebox{1}{
		\includegraphics[width=0.08\textwidth]{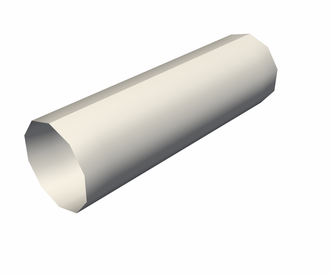}
		\includegraphics[width=0.08\textwidth]{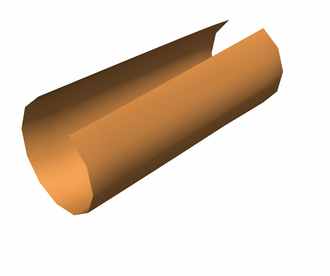}
		\includegraphics[width=0.08\textwidth]{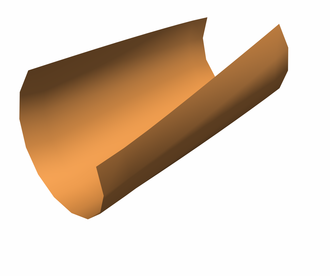}
		\includegraphics[width=0.08\textwidth]{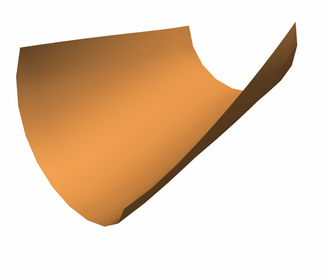}
		\includegraphics[width=0.08\textwidth]{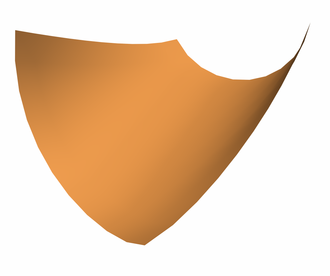}
		\includegraphics[width=0.08\textwidth]{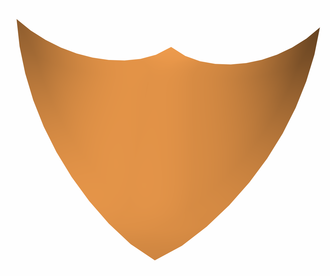}
		\includegraphics[width=0.08\textwidth]{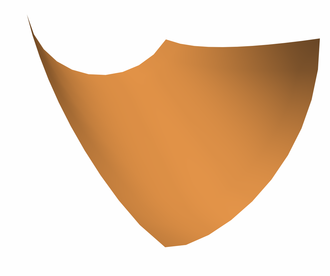}
		\includegraphics[width=0.08\textwidth]{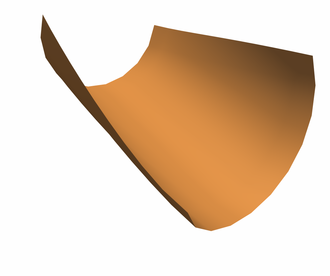}
		\includegraphics[width=0.08\textwidth]{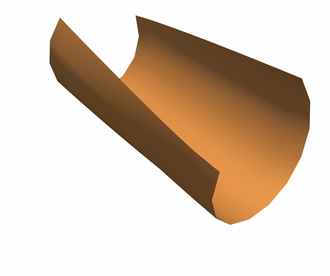}
		\includegraphics[width=0.08\textwidth]{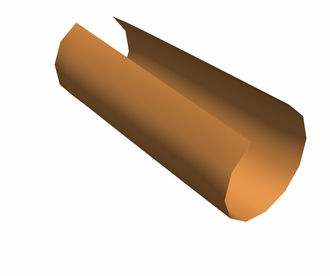}
		\includegraphics[width=0.08\textwidth]{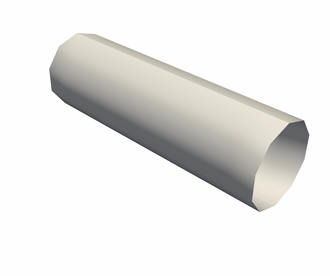}
	}\\[0.5em]
	\scalebox{1}{
		\begin{tikzpicture} 
		\begin{axis}[
		ymajorticks=false,
		axis y line=none,
		axis x line*=bottom,
		enlargelimits=0.05,
		enlarge x limits=0.001,
		ybar interval=0.8, 
		ymin=0.0015,
		ymax=0.0153,
		width=0.95\textwidth,
		height=2cm 
		]
		\addplot[fill=myOrange!80!gray] coordinates {(0,0.003534) (1,0.003535) (2,0.003537) (3,0.003549) (4,0.003602) (5,0.003607) (6,0.003551) (7,0.003538) (8,0.003535) (9,0.003535) (10,0.003535)};
		\end{axis} 
		\end{tikzpicture}
	}\\[0.5em]
	\includegraphics[width=\columnwidth]{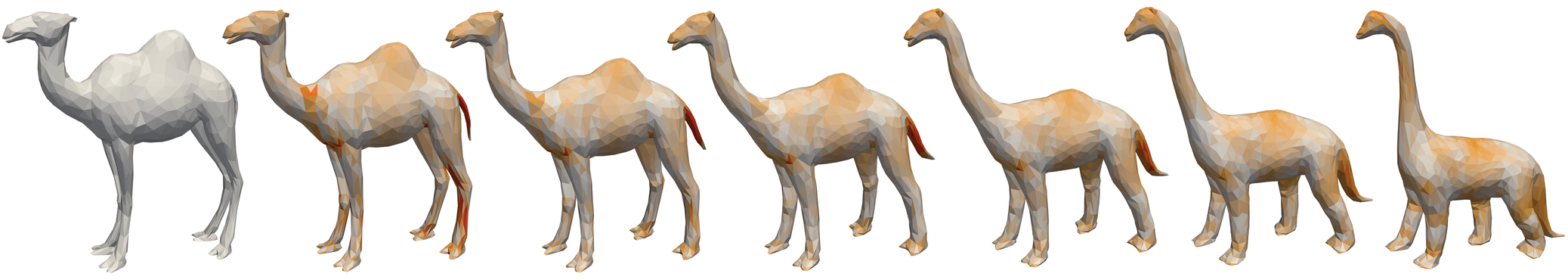}\\[0.5em]
	\scalebox{1}{
		\hspace{1em}
		\begin{tikzpicture} 
		\begin{axis}[
		ymajorticks=false,
		axis y line=none,
		axis x line*=bottom,
		enlargelimits=0.05,
		enlarge x limits=0.001,
		ybar interval=0.8, 
		ymin=0.0015,
		ymax=0.0103,
		width=0.9\textwidth,
		height=2cm 
		]
		\addplot[fill=myOrange!80!gray] coordinates {(0,0.003534) (1,0.003535) (2,0.003537) (3,0.003549) (4,0.003602) (5,0.003607) (6,0.003551)};
		\end{axis} 
		\end{tikzpicture}
	}
	\caption{Bottom: discrete geodesic in NRIC with input data from \cite{AmRo18} and membrane distortion as colormap ($0\hspace{1mm}$\protect\resizebox{.1\linewidth}{!}{\protect\includegraphics{paperfolds/colorbar_red_log.png}}$\hspace{1mm}\geq1$); 
		above: Linear interpolation in ambient space $\R^{2\numE}$ as in \cite{FrBo11} with energy distribution (green) vs. geodesic interpolation on $\manifold$ with constant energy distribution (orange).} 
	\label{fig:geodesicExample}
\end{figure}

The concept of discrete geodesics directly translates to the NRIC manifold and the path energy in \eqref{eq:discretePathEnergy} can be considered as an objective functional in \autoref{eq:genericOpt}. 
Note, however, that this increases the number of free variables substantially. 
In \autoref{fig:geodesicExample}, we show different time-discrete geodesics in NRIC where we use the quadratic deformation energy \autoref{eq:quadEnergy} in \autoref{eq:discretePathEnergy}. 
In particular, we compare for end shapes being two oppositely bent plates our NRIC geodesic (orange) to the linear interpolation (green) in the embedding space $\R^{2\numE}$, which corresponds to a naive transfer of the projection approach by \citet{FrBo11}. 
As indicated by the histogram plots, the discrete constant speed property can only be obtained for the NRIC formulation.

\begin{figure}[ht]
	\centering
	\includegraphics[width=\columnwidth]{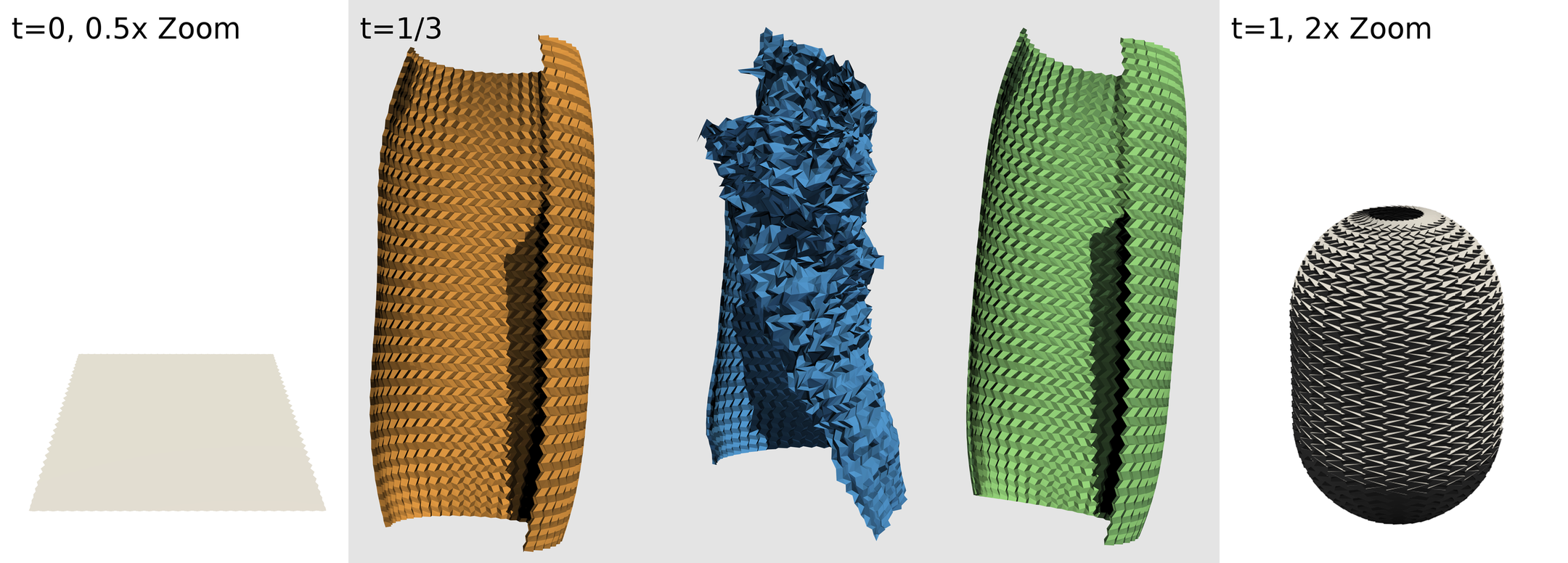}
	\caption{Intermediate shapes at $t=1/3$ of a discrete geodesic between two perfectly isometric end shapes (grey) taken from \citet{DuVoTaMa16} obtained via NRIC optimization (orange) and vertex-based methods as in \cite{HeRuWa12} (blue) resp. \cite{FrBo11} (green). 
		Note that we preserve the isometry due to our hard length constraints. In contrast, the vertex-based methods get either stuck in local minima (blue) or reveal artifacts such as unnatural asymmetries (green).}
	\label{fig:pillGeodesic}
\end{figure}

\begin{figure*}[ht]
	\includegraphics[width=0.7\linewidth]{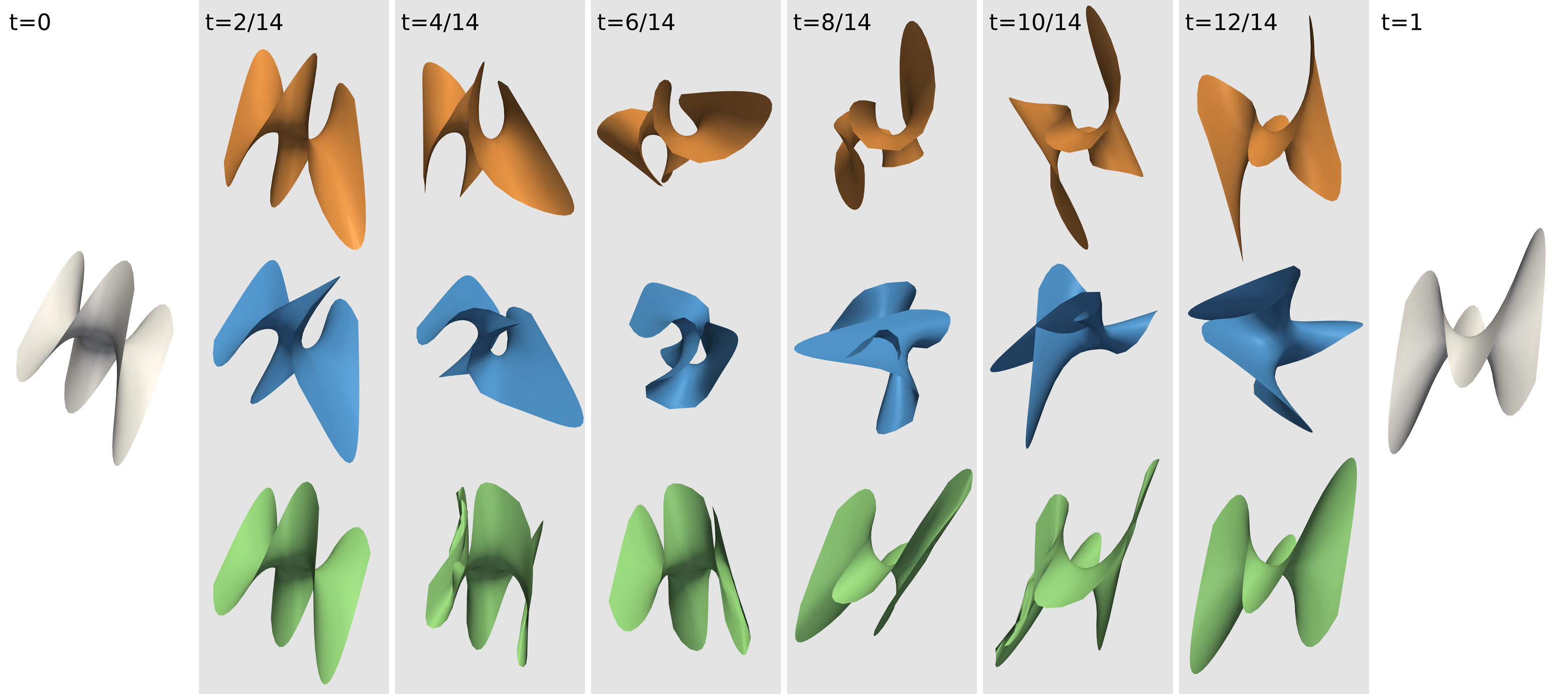}
	\begin{tikzpicture} 
	\begin{axis}[ticks=none,
	xlabel={\(t\)},
	axis y line=right,
	axis x line=bottom,
	enlarge x limits={abs=15pt},
	x axis line style=-,
	ylabel={\( \lVert l(t) - l^\ast \rVert_2 \)}, 
	y label style={yshift=1.5em},
	x tick style={below},
	grid=none,
	ymin=-0.005,
	ymax=0.4,
	ybar stacked,
	bar width=13pt,
	width=0.3\linewidth,
	height=5cm,
	compat=1.9
	]
	\addplot +[forget plot,fill=myOrange,color=myOrange] coordinates {(2,-1) (4,-1) (6,-1) (8,-1) (10,-1) (12,-1) };
	\addplot[fill=myOrange,color=myOrange] coordinates {(2,1) (4,1) (6,1) (8,1) (10,1) (12,1) };
	
	\addplot[fill=myBlue!70!gray,color=myBlue!70!gray] coordinates {(2,0.093645) (4,0.147862) (6,0.166342) (8,0.15664) (10,0.133562) (12,0.0788468)};
	\addplot[fill=myGreen!70!gray,color=myGreen!70!gray] coordinates {(2,0.0148586) (4,0.198558) (6,0.133258) (8,0.133403) (10,0.198562) (12,0.0148603) };
	
	%
	
	\end{axis} 
	\end{tikzpicture}
	\caption{Isometric geodesic paths. Left input shape $(\len^\ast, \dih^\ast)$ as hyperbolic monkey saddle, right input shape is the reflection $(\len^\ast, -\dih^\ast)$. 
		Comparison of discrete geodesics computed in NRIC (orange, perfectly isometric) and by methods based on nodal positions (\ie \citet{HeRuWa12} (blue) and \citet{FrBo11}). 
		The latter approaches are not able to resolve pure isometric geodesics as indicated by a histogram of varying lengths on the right. See also supplementary video.}
	\label{fig:monkeySaddleGeodesic}
\end{figure*}

Furthermore, we can combine the computation of time-discrete geodesics with further constraints on the coordinates, \eg to simulate isometric deformation paths. 
For example, in \autoref{fig:pillGeodesic} we compare the computation of (almost) \emph{isometric} geodesic paths between perfectly isometric end shapes taken from \citet{DuVoTaMa16} to vertex-based methods. 
A similar example is shown in \autoref{fig:monkeySaddleGeodesic}, where the first input shape $\z_0 = (\len^\ast, \dih^\ast)$ describes a hyperbolic monkey saddle and the second input shape is given by a reflection $\z_K = (\len^\ast, -\dih^\ast)$ of the saddle. 
\autoref{fig:monkeySaddleGeodesic} demonstrates that our approach is able to realize a perfectly isometric deformation path (orange) by enforcing $\len_k = \len^\ast$ for all $0<k<K$, whereas vertex based optimization methods fail. 

\paragraph*{Timings}
Lastly, let us discuss the runtimes of the proposed method, where all computations were performed using a desktop computer with an Intel(R) Core(TM) i7-4790 CPU and 16 GB RAM.
In our framework, we use the Eigen library \cite{Eigen} for linear algebra tasks and CHOLMOD \cite{CHOLMOD} for the Cholesky decomposition.
At first, we list timings for the reconstruction performed without parallelization. 
As a representative example, we report on timings measured on the Dyna dataset \cite{PoRoMaBl15} (\cf \autoref{fig:dynaReconstruction}) where $\numV \approx 6.9k$. 
The generation of (MST) or (SPT) takes about 12ms, the generation of a spanning tree via (BFS) takes about 1ms. 
The traversal of a spanning tree takes about 5ms, and one Gau\ss-Newton iteration is done in 330ms (with 80\% spending in the linear solver).
Next, computing the entries of the Hessian of the constraint functional $\qint$ requires, again without parallelization, 16ms for the discrete surfaces considered in \autoref{fig:paperFoldingCompareGauss} and 112ms for the Dyna dataset considered in \autoref{fig:strangulation}.
Detailed timings for the optimization described in \autoref{sec:varProblems} are listed in \autoref{table:runtimes}.
Note, that the evaluation of the augmented Lagrangian $\lagrange$ and its derivatives requires substantially more time than computing the entries of $D^2_\z \qint$. 
This originates from computing the square of $D_\z \qint$ and assembling the Hessian in CSR format because these operations do not benefit from parallelization in our current implementation.

Compared to computations in nodal positions, our method requires more memory due to the increased number of primal degrees of freedom.
However, because this number is approximately twice the number of nodal positions the total memory consumption only increases by a constant factor of approximately four.

\begin{table}[t]
	\centering	
	\renewcommand{\arraystretch}{1.1}	
	\begin{tabular}{|lc||ccccc|}
		\hline
		\multicolumn{2}{|c||}{example} & \multicolumn{2}{c}{iterations} & \multicolumn{3}{c|}{avg.\ times per Newton iteration} \\ 
		Figure & \(N\)  & aug.\ Lagrange & Newton & evaluation & solve & line search \\ 
		\hline \hline
		\ref{fig:strangulation} (left, avg.) & 5220  & 18 & 83 & 30 ms & 23 ms & 3 ms \\ 
		\ref{fig:strangulation} (right) & 41328 & 14 & 121 & 281 ms & 201 ms & 15 ms \\ 
		\ref{fig:handsAverage} (avg.) & 36552 & 13 & 203 & 263 ms & 244 ms & 18 ms \\ 
		\ref{fig:paperFoldingCompareGauss} & 6272  & 83 & 477 & 34 ms & 34 ms & 7 ms \\ 
		\ref{fig:paperFoldingCompare} (right) & 24832 & 46 & 454 & 173 ms & 84 ms & 7 ms \\ 
		\ref{fig:geodesicExample} (top) & 5760 & 10 & 15 & 16 ms & 69 ms & 2 ms \\ 
		\ref{fig:geodesicExample} (bottom) & 36000 & 12 & 65 & 338 ms & 732 ms & 10 ms \\ 
		\ref{fig:pillGeodesic} & 69936 & 11 & 44 & 626 ms & 301 ms & 22 ms \\ 
		\ref{fig:monkeySaddleGeodesic} & 32240 & 11 & 173 & 218 ms & 193 ms & 5 ms \\ 
		\hline
	\end{tabular}
	\caption{
		Performance statistics of our approach on the different examples shown before. 
		From left to right: number of degrees of freedom, number of iterations of the augmented Lagrange method, total number of Newton iterations, average time for evaluation of function and derivatives per Newton iteration, average time for computing $\shift$ and solving the linear system per Newton iteration, and average time for line search per Newton iteration.
	}
	\label{table:runtimes}
\end{table}

\section{Conclusion}\label{sec:discussion}
We introduce a framework that allows us to pose and solve geometric optimization problems in terms of NRIC. 
The framework is built on several novel concepts. 
First, we introduce a Riemannian structure for the NRIC manifold stemming from an implicit description via integrability conditions and a physically-motivated nonlinear elastic energy. 
In particular, we demonstrate how the notion of a tangent space can be used to identify infinitesimal isometric variations. 
Second, we present an approach based on the augmented Lagrange method and a modified line search for solving generic optimization problems in NRIC. 
Third, we develop a hybrid algorithm for the reconstruction of nodal positions from length and angle coordinates that uses a mesh traversal to initialize a Gau\ss--Newton solver. 
We tested our framework on different problems including shape interpolation and paper folding. 
A particular strength is the simulation of true isometric deformations---a task where well-established vertex-based methods often fail. 

\paragraph*{Limitations and challenges} 
We see great potential in using NRIC for geometric optimization problems and expect that the techniques we present will be further developed. 
We plan to formulate an extended geodesic shape space calculus (\cf \cite{HeRuSc14}) including geodesic extrapolation and parallel transport in NRIC and expect to profit from the rigid motion invariance of the coordinates and their robustness for near-isometric deformation.  
In the context of a statistical analysis of shapes, our NRIC formulation enables direct processing of input data without an a priori rigid co-registration. 
To this end, our NRIC manifold is a natural starting point for the development of a corresponding Riemannian principal component analysis. 

Though our experiments demonstrate the benefits of NRIC-based optimization, our current framework has several limitations and poses challenges in making the optimization more efficient. 
First, the current implementation can only handle simply connected surfaces. 
An extension to higher-genus surfaces would require to include integrability conditions along non-contractible paths that generate the fundamental group. 
This would lead to more global constraints in our optimization problems.
Typical examples of surfaces in geometric modelling have only a small number of generators of the homology group. 
However, this necessity of complicated constraints is a general limitation of our method compared to nodal positions.

Second, a fundamental challenge is to reduce the number of degrees of freedom and integrability conditions. 
Our current framework works with $2\numE$ variables and $3\numV$ integrability conditions per shape. 
This implies a larger number of variables compared to optimization in nodal positions, which in turn means increased memory requirement and more costly iterations. 
Here, it might be worthwhile to explore model reduction approaches.
Furthermore, the triangle inequality constraints are in general challenging to take care of in the implementation.
We found in all our experiments that the proposed adapted line search, especially in conjunction with the nonlinear deformation energy, was able to handle them robustly.

Finally, we aim to account for point constraints in our NRIC-based optimization. 
These type of constraints frequently appear, for instance, in physical simulations as forces or boundary conditions.
This could be accomplished by performing a partial reconstruction of the points with attached constraints using an explicit formula that results from tracing the paths in \autoref{alg:frame_direct}.
Then the derivatives of the explicit formula need to be computed with respect to NRIC to enable their use in optimization problems which might be a feasible task for modern automatic differentiation frameworks.
Nonetheless, this would introduce highly nonlinear and nonlocal terms to the optimization potentially limiting the performance of our method.
This introduces the challenge of devising different ways to combine NRIC-based modeling with point constraints.

\appendix

\section{Rotations $q_{ij}$}
As noted above, each induced transition rotation depends on one dihedral angle and on the three edge lengths of a triangle.
We simplify our notation for the computation of the partial derivatives and define
\begin{equation}
\hat q(\theta,a,b,c) := q_0(\theta)\, q_2\left(\arccos\left(\frac{a^2+b^2-c^2}{2 ab}\right)\right),
\end{equation}
where \(a,b,\) and \(c\) are the edge lengths and \(\theta\) is the dihedral angle.
To simplify the notation even further in the following, we define the rational function \(Q(a,b,c) = \frac{a^2+b^2-c^2}{2 ab}\).
The angle of the rotation around the second standard basis vector is given by the law of cosines and using trigonometric formulas we can thus simplify its matrix representation to
\begin{equation*}
q_2\left(\arccos\, Q\right) = \sqrt{\frac{1+Q}{2}}+ \sqrt{\frac{1-Q}{2}}\, \boldsymbol{k} .
\end{equation*}
If we multiply this with the rotation around the zeroth standard basis vector we arrive at
\begin{equation*}
\hat q(\theta,a,b,c) = \cos \frac{\theta}{2} \sqrt{\frac{1+Q}{2}} +   \sin \frac{\theta}{2} \sqrt{\frac{1+Q}{2}}\, \boldsymbol{i} - \sin \frac{\theta}{2} \sqrt{\frac{1-Q}{2}}\, \boldsymbol{j} + \cos \frac{\theta}{2} \sqrt{\frac{1-Q}{2}} \, \boldsymbol{k} .
\end{equation*}
With this representation at hand, it is now possible to compute its first and second derivative.
In particular, this is a viable task for a symbolic differentiation tool and we refer to the supplementary material for the results.
\section{Local membrane energy}
To understand the nonlinear membrane energy $\W_\mem$ better and prepare the computation of its derivatives, we study in this section the contribution of a single triangle $\face \in \faces$.
To this end, let $a,b,c$ be the edge lengths of $\face$ in the undeformed configuration $\z$ and $\tilde{a}, \tilde{b}, \tilde{c}$ the corresponding edge lengths of the deformed $\tilde \z$.
Our goal now is to express $\area_\face \cdot W_\mem(\mathcal G[\z, \tilde \z]|_\face)$ in terms of these lengths.

We start with the components of $\dfFF$  as given in \autoref{eq:discFirstFundForm}.
The diagonal entries are of course simply given by $b^2$ and $a^2$.
For the off-diagonal entries, recall 
\[\langle E_1(f), E_2(f)\rangle = \lVert  E_1(f) \rVert \, \lVert  E_2(f) \rVert \cos \left(\angle (E_1(f), E_2(f))\right), \] 
where by the law of cosines we have $\cos \angle (E_1(f), E_2(f)) = \frac{a^2+b^2-c^2}{2ab}$.
Hence, we get the representation 
\begin{equation}
	\hat \dfFF(a,b,c) = \begin{pmatrix}
	b^2 & -\frac{1}{2} (a^2+b^2-c^2) \\
	-\frac{1}{2} (a^2+b^2-c^2) & a^2
	\end{pmatrix},
\end{equation}
where we use the hat to indicate the local representation of an object.
From this, we immediately compute the determinant as
\begin{equation}
	\det \hat \dfFF(a,b,c) = \frac{1}{4} (a+b+c)(-a+b+c)(a-b+c)(a+b-c) = 4\, \area_\face(a,b,c)^2,
\end{equation}
with $\area_\face(a,b,c)$ the triangle area obtained via Heron's formula.
The local representation of the discrete distortion tensor is given by
\begin{equation}
	\hat \dDisTen (a,b,c,\tilde{a},\tilde{b},\tilde{c}) = \hat\dfFF(a,b,c)^{-1} \hat\dfFF(\tilde{a}, \tilde{b}, \tilde{c}).
\end{equation}
By the multiplicativity of the determinant, we obtain 
\begin{equation}
	\det \hat \dDisTen = (\det \hat \dfFF(a,b,c))^{-1} \hat \dfFF(\tilde{a}, \tilde{b}, \tilde{c})  = \frac{\area_\face(\tilde{a}, \tilde{b}, \tilde{c})^2}{\area_\face(a,b,c)^2}
\end{equation}
for the determinant of the distortion tensor. 
Computing the trace requires in contrast an explicit representation of $\hat \dDisTen$ and we finally obtain
\begin{equation}
	\tr \hat \dDisTen = \frac{1}{8\, \area_\face(a,b,c)^2} \left( \tilde{a}^2(-a^2+b^2+c^2)  +  \tilde{b}^2 (a^2-b^2+c^2) + \tilde{c}^2 (a^2+b^2-c^2) \right).
\end{equation}
Together, we have assembled all components to write the contribution of $\face$ as 
\begin{equation}
	\area_\face(a,b,c) \cdot W_\mem(\hat \dDisTen (a,b,c,\tilde{a},\tilde{b},\tilde{c}))
\end{equation}
with $W_\mem$ the energy density from \autoref{def:ltheta_membrane_energy}.
We provide the derivatives of the membrane energy along with the derivatives of the bending energy in the supplementary material.

\section{Optimization algorithms}
In this section, we will discuss the algorithmic details of the optimization procedure introduced in \autoref{sec:varProblems} for our generic problem \autoref{eq:genericOpt}.
First, we start with detailed description of the augmented Lagrange method adapted from \cite{NoWr06} in \autoref{alg:augmentedLagrange}, which provides all parameters related to the increase of the penalty parameter and the update of the Lagrange multiplier estimates.
In almost all examples, our default parameters $\pen^0 = 10,\ \mult^0 = 0,\ \pen_{+} = 100,$ and $\eta_{+}=0.9$ worked reasonably well.
Only in the paper folding examples (\cf \autoref{fig:paperFoldingCompareGauss} and \autoref{fig:paperFoldingCompare}), we decreased $\pen_{+}$ and $\eta_{+}$ because we noticed this leads to local minima with lower energy values as was discussed in \autoref{sec:results}.
\begin{algorithm}[h]
	\caption{Augmented Lagrange method for \autoref{eq:genericOpt}, \cite[Alg.\ 17.4]{NoWr06}}
	\label{alg:augmentedLagrange}
	\begin{algorithmic}[1]
		\Require Initial NRIC $\z^0$, initial penalty $\pen^0$, initial multipliers $\mult^0$, constraint tolerance $\varepsilon_\qint$, optimality tolerance $\varepsilon_L$, penalty increase factor $\pen_{+}$, tolerance increase exponent $\eta_{+}$ 
		\Ensure Approximate solution $\z^k$ to \autoref{eq:genericOpt}
		\ForAll{$k=0,1,2,\ldots,k_{max}$}
		\State Compute approximate minimizer $\z^{k+1}$ of $\lagrange(\,\cdot\,, \mult^k, \pen^k)$ with $\lVert D_z L(\z^{k+1}, \mult^k, \pen^k)\rVert_2 \leq \omega_k$
		\If{$\lVert \qint(\z^{k+1}) \rVert_\infty \leq \varepsilon_\qint$ and $\lVert D_z L(\z^{k+1}, \mult^k, \pen^k)\rVert_2 \leq \varepsilon_L$} 
			\Comment{Stopping criterion}
			\State \Return $\z^{k+1}$
		\EndIf
		\If{$\lVert \qint(\z^{k+1}) \rVert_\infty \leq \eta^k$} \Comment{Recompute multiplier}
			\State $\mult^{k+1} = \mult^k - \pen^k \qint(z^{k+1})$ 
			\State $\pen^{k+1} = \pen^k$
			\State $\eta^{k+1} = \max ({\eta^k} / {(\pen^{k+1})^{\eta_+}}, \varepsilon_\qint)$
			\State $\omega^{k+1} = {\omega^k} / {{\pen^{k+1}}}$
		\Else \Comment{Increase penalty parameter}
			\State $\mult^{k+1} = \mult^k$
			\State $\pen^{k+1} = \pen_+\cdot \pen^k$
			\State $\eta^{k+1} = \max (1 / {(\pen^{k+1})^{0.1}}, \varepsilon_\qint)$
			\State $\omega^{k+1} = 1 / {{\pen^{k+1}}}$
		\EndIf
		\EndFor
		
	\end{algorithmic}
\end{algorithm}

To compute the approximate minimizer, we use the variation of the Newton-type method introduced before which requires the first and second derivatives of the augmented Lagrangian $\lagrange$.
The first derivative of $\lagrange$ is given by
\begin{equation*}
	D_\z \lagrange(\z, \mult, \pen) = D_\z \energy(\z) - D_\z \qint(\z) \cdot \mult + \pen\, D_\z \qint(z)^T  \qint(z) 
\end{equation*}
and hence the second derivative turns out to be
\begin{equation*}
	D^2_\z\lagrange(\z, \mult, \pen) = D^2_\z \energy(\z) + D^2_\z \qint(\z) \cdot \left(\pen \, \qint(z) -\mult \right) + \pen\, D_\z \qint(z)^T  D_\z \qint(z),
\end{equation*}
where 
\begin{equation*}
	D_\z^2 \qint \cdot \left(\pen \, \qint(z) -\mult \right) = \left(\sum_{\vertex \in \vertices_0}  \partial_{\z_l} \partial_{\z_k} \qint_\vertex \cdot \left(\pen \, \qint_\vertex(z) -\mult_\vertex \right) \right)_{l,k = 1,\ldots, 2\numE}\,.
\end{equation*}
Again, the detailed derivatives of the integrability constraints and of the energy are provided in the supplementary material.
Now, we can provide all the steps of this method in an integrated fashion as \autoref{alg:newton}.
Again, in almost all examples our default parameters $\shift_{+} = 10$ and either $\beta=10^{-3}$ or $\beta=10^{-4}$ worked well.
\begin{algorithm}[h]
	\caption{Line search Newton-type method for $L$, \cite[Alg.\ 3.1, 3.2, 3.3]{NoWr06} }
	\label{alg:newton}
	\begin{algorithmic}[1]
		\Require Initial NRIC $\z^k$, gradient tolerance $\omega^k$, initial/minimal shift $\beta$, shift increase factor $\shift_+$
		\Ensure Approximate minimizer $\z^{k+1}$ of $\lagrange(\, \cdot \, , \mult^k, \pen^k)$
		\State Set $\z^{k,0} = \z^k$
		\For{$j=0,1,2,\ldots,j_{max}$}
		
		\State Evaluate $\lagrange(\z^{k,j}, \mult^k, \pen^k),\ D_\z \lagrange(\z^{k,j}, \mult^k, \pen^k),$ and $D^2_\z \lagrange(\z^{k,j}, \mult^k, \pen^k)$ \Comment{Evaluation}
		
		\If{$\lVert D_\z \lagrange(\z^{k,j}, \mult^k, \pen^k) \rVert \leq \omega^k$ } \Comment{Stopping criterion}
			\State \Return $\z^{k+1} = \z^{k,j}$
		\EndIf
		
		\If{$\min_i \left(D^2_\z \lagrange(\z^{k,j}, \mult^k, \pen^k)\right)_{ii} > 0$} \Comment{Determine shift}
			\State $\shift_j = 0$
		\Else
			\State $\shift_j = -\min_i \left(D^2_\z \lagrange(\z^{k,j}, \mult^k, \pen^k)\right)_{ii} + \beta$
		\EndIf
		
		\Loop
			\State Attempt Cholesky decomposition of $D^2_\z \lagrange(\z^{k,j}, \mult^k, \pen^k) + \shift_j \Id$
			\If{factorization succeeds}
				\State \textbf{stop}
			\Else
				\State $\shift_j \leftarrow \max(\shift_+ \cdot \shift_j,\ \beta)$
			\EndIf
		\EndLoop
		
		\State Solve $\left(D^2_\z \lagrange(\z^{k,j}, \mult^k, \pen^k) + \shift_i \Id\right) d_j = -D_\z \lagrange(\z^{k,j}, \mult^k, \pen^k)$ \Comment{Descent direction}
		
		\Repeat  \Comment{Line search}
			\State $\alpha_j = 0.5 \alpha_j $
		\Until{$\lagrange(\z^{k,j}+\alpha_j d_j, \mult^k, \pen^k) \leq \lagrange(\z^{k,j}+\alpha_j d_j, \mult^k, \pen^k) + 0.1 \alpha_j \, D_\z \lagrange(\z^{k,j}, \mult^k, \pen^k)^T d_j $}
		
		\State Set $\z^{k,j+1} = \z^{k,j}+\alpha_j d_j$
		
		\EndFor
	\end{algorithmic}
\end{algorithm}
\section{Direct reconstruction algorithm}
In Section \ref{sec:recon}, we have already outlined how the frames can be constructed iteratively using the transition rotations induced by $\z\in\manifold$.
To complete the description of the reconstruction algorithm, we also need to detail how to construct the nodal positions.
To this end, consider a face with a given discrete frame \(\frame\) and target edge lengths $\len_0, \len_1, \len_2$.
Then we obtain embedded edge vectors $E_i \in \R^3$ with $\|E_i\| = \len_i$ for $i = 0,1,2$ by
\begin{equation}
\label{eq:edgeVectors}
\begin{gathered}
E_1 = \len_0 \frame \begin{pmatrix}1\\0\\0\end{pmatrix}\!, \ 
E_2 = \len_1 \frame \begin{pmatrix}-\cos \gamma_3\\\sin \gamma_3\\0\end{pmatrix}\!,
E_3 = \len_2 \frame \begin{pmatrix}-\cos \gamma_2\\ -\sin \gamma_2\\0\end{pmatrix}\!,
\end{gathered}
\end{equation}
and finally nodal positions \((\pos_i)_{i=0,1,2}\) such that \(E_i = \pos_{i-1} - \pos_{i+1}\). 
Note that the inner angles $\gamma_2$ and $\gamma_3$ can be obtained from the edge lengths using the law of cosines.
The complete reconstruction algorithm is summarized in Algorithm \ref{alg:frame_direct}. 
Note that in practice, we need to construct at most one nodal position per face (except for the first face $\face_0$), and often even none as the positions are already determined.
\begin{algorithm}[ht]
	\caption{Direct frame-based reconstruction}
	\label{alg:frame_direct}
	\begin{algorithmic}[1]
		\Require $\connec = (\vertices,\edges)$, $\z\in\R^{2\numE}$,  \(\pos_0 \in \R^3\) and pair $(\frame_0, \face_0)$
		\Ensure Nodal positions \(\pos\) with $\projZ(\pos) \approx \z$ 
		\State Evaluate the discrete integrability map \(\dint_\vertex(\z)\) for all \(\vertex \in \vertices\)
		\State Define edge weights \eqref{eq:edgeWeights} for \(e=(vv') \in \edges\)
		\Statex  \hspace{-\algorithmicindent}  \textbf{Alternative MST:}
		\setcounter{ALG@line}{2}	
		\State Construct minimal spanning tree of the dual graph based on \(w\)
		\Statex \hspace{-\algorithmicindent} \textbf{Alternative SPT:}
		\setcounter{ALG@line}{2}
		\State Construct shortest path tree of the dual graph based on \(w\)
		\State Traverse the dual graph of $\connec$ following the constructed tree
		\ForAll{dual edges \((\face_i, \face_j)\)}
		\State Construct induced transition rotation \(R_{ij}\)
		\State Construct frame \(\frame_j = \frame_{i} R_{ij}\)
		\State Compute interior angles \(\gamma_1, \gamma_2, \gamma_3\) of \(\face_j\)
		\State Determine edge vectors \(E_1, E_2, E_3\) according to \eqref{eq:edgeVectors}
		\State Determine nodal positions for \(\face_j\) 
		\EndFor
	\end{algorithmic}
\end{algorithm}


\section*{Acknowledgements}
The authors thank Laszlo Bardos from \href{http://www.cutoutfoldup.com}{cutoutfoldup.com} for the photo of Steffen's polyhedron used in \autoref{fig:steffens}, furthermore the authors of \cite{BoVoGoWa16} and the American Society of Mechanical Engineers for the photo in \autoref{fig:origami}.
Moreover, we thank Etienne Vouga and Levi Dudte for the meshes used in \autoref{fig:pillGeodesic}, Friedrich Bös for information on the construction of the Origami cylinder in \autoref{fig:origami}, and Carlos Rojas for the animal meshes used in \autoref{fig:geodesicExample}.
The authors from Bonn gratefully acknowledge the support provided by the Austrian Science Fund (FWF) through project NFN S117 and by the Deutsche Forschungsgemeinschaft (DFG) through the Hausdorff Center for Mathematics (GZ 204711, Project ID 390685813).


\bibliographystyle{elsarticle-num-names} 
\bibliography{sgp2019references}

\end{document}